\documentclass[11pt,a4paper]{article}
\usepackage{amssymb,amsmath, amsfonts}
\usepackage{graphicx,graphics}
\usepackage{mathtools}
\usepackage{bbold}
\usepackage[english]{babel}
\usepackage[utf8]{inputenc}
\usepackage{epsfig,url}
\usepackage{bbm,theorem}
\usepackage{a4wide}
\usepackage{color}
\usepackage{enumerate}
\usepackage{calrsfs}
\usepackage[bookmarks=true]{hyperref}
\usepackage{bookmark}
\usepackage{verbatim}
\DeclareMathAlphabet{\pazocal}{OMS}{zplm}{m}{n}
\newtheorem{theorem}{Theorem}[section]
\newtheorem{definition}[theorem]{Definition}

\newtheorem{lemma}[theorem]{Lemma}
\newtheorem{proposition}[theorem]{Proposition}

\newtheorem{remark}[theorem]{Remark}

\numberwithin{equation}{section}
\numberwithin{theorem}{section}
\newcommand{\qed}{\hfill$\Box$}
\newenvironment{proof}{\begin{trivlist}\item[]{\em Proof:}\/}{\qed\end{trivlist}}
\newenvironment{proofof}[1][Proof]{\noindent \textit{#1.} }{\ \qed}
\newcommand{\R}{{\mathbb R}}
\newcommand{\K}{{\mathbb K}}
\DeclareMathOperator*{\supp}{supp}

\begin{document}

\title{Asymptotic localization in multicomponent mass conserving  coagulation equations \\ \ \\ 
\large {\bf Localization in multicomponent coagulation}}

\author{Marina A. Ferreira, Jani Lukkarinen, Alessia Nota, Juan J. L. Vel\'azquez}
\maketitle

\begin{abstract}
In this paper we prove that the time dependent solutions of a large class of
 Smoluchowski coagulation equations for multicomponent systems concentrate
along a particular direction of the space of cluster   compositions for long times. The
direction of concentration is determined by the initial distribution of
clusters. These results allow to prove the uniqueness and global stability of
the self-similar profile  with finite mass in the case of coagulation
kernels which are not identically constant, but are constant along any
direction of the space of cluster compositions.

\end{abstract}

\bigskip

\textbf{Keywords:}  multicomponent Smoluchowski's equation;  localization;  time-dependent solutions; self-similarity; stability.

\tableofcontents

\section{Introduction}

\subsection{Motivation}

In this paper we are concerned with two classes of multicomponent Smoluchowski
coagulation equations, namely the discrete equation and continuous equation.
In the discrete case, clusters are characterized by the composition vector
$\alpha=\left(  \alpha_{1},\alpha_{2},\ldots,\alpha_{d}\right)  \in
\mathbb{N}_{0}^{d}\backslash\{ O\}$ consisting of $d$ different monomer types,
where $\mathbb{N}_{0}=\left\{  0,1,2,3,\ldots\right\}  $ and $O=\left(
0,0,\ldots,0\right)  .$ The concentration $n_{\alpha}(t)$ of particles of
composition $\alpha$ at time $t\geq 0$ is governed by the following equation
\begin{equation}
\partial_{t}n_{\alpha}= \kappa_d[n_\alpha],\quad\quad \kappa_d[n_\alpha] := \frac{1}{2}\sum_{\beta<\alpha}K_{\alpha-\beta,\beta
}n_{\alpha-\beta}n_{\beta}-n_{\alpha}\sum_{\beta>O}K_{\alpha,\beta}n_{\beta}\,.
\label{B1}
\end{equation}
Given $\alpha=\left(  \alpha_{1},\alpha_{2},\ldots,\alpha_{d}\right)  $ and
$\beta=\left(  \beta_{1},\beta_{2},\ldots,\beta_{d}\right)  $ we write
$\beta<\alpha$ to indicate that $\beta_{k}\leq\alpha_{k}$ for all
$k=1,2,\ldots,d$, and in addition $\alpha\neq\beta.$ The collision kernel
$K_{\alpha,\beta},$ which we assume to satisfy the symmetry condition $K_{\alpha,\beta
}=K_{\beta,\alpha},$ describes the coagulation rate between clusters with
compositions $\alpha$ and $\beta.$ Equation \eqref{B1} was first proposed in
\cite{S16} in the case of particles described by a single component,
corresponding here to the case $d=1$.

In the following, we will most often work with the composition spaces ${\mathbb{R}}_{\ast}=(0,\infty)$ and ${\mathbb{R}}_{\ast}%
^{d}=[0,\infty)^{d}\backslash\{O\}$, with $O=\left(  0,0,...,0\right)$.
The continuous version of  equation (\ref{B1}) is then given by
\begin{equation}
\partial_{t}f\left(  x,t\right)  = \K_d[f](x,t)  \,,\ \ x\in\mathbb{R}_{\ast}%
^{d}, \quad t\geq 0
\label{B4}
\end{equation}
where 
\begin{equation}
 \K_d[f](x,t)  : = \frac{1}{2}\int_{\left\{  0<\xi<x\right\}
}d\xi K\left(  x-\xi,\xi\right)  f\left(  x-\xi,t\right)  f\left(
\xi,t\right)  -\int_{{\mathbb{R}}_{\ast}^{d}}d\xi K\left(  x,\xi\right)
f\left(  x,t\right)  f\left(  \xi,t\right)\label{coagKern}
\end{equation}
and for $x=\left(  x_{1},x_{2},\ldots,x_{d}\right)  ,$ $y=\left(
y_{1},y_{2},\ldots,y_{d}\right)  $ we use the previously introduced comparison
notation: $x<y$ whenever $x\leq y$ componentwise and $x\neq y$. In
particular, we thus have
\[
\int_{\left\{  0<\xi<x\right\}  }d\xi=\int_{0}^{x_{1}}d\xi_{1}\int_{0}^{x_{2}%
}d\xi_{2}\cdots\int_{0}^{x_{d}}d\xi_{d}\,.
\]
We will also 
assume that $
K(x,y)=K(y,x),\ K(x,y)\geq0.
$

Notice that (\ref{B1}) can be considered a particular case of (\ref{B4}) for
measure solutions of (\ref{B4}) with the form%
\begin{equation}
f\left(  x,t\right)  =\sum_{\alpha\in\mathbb{N}_{0}^{d}\setminus \{ O\}}n_{\alpha
}\delta\left(  x-\alpha\right).  \label{B5}%
\end{equation}
In this case, $f$ in \eqref{B5} solves \eqref{B4} if $(n_{\alpha})_{\alpha}$
solves \eqref{B1} with the kernel $K_{\alpha,\beta}=K\left(  \alpha
,\beta\right)$. We have chosen to state the results separately for each
equation because the results for the discrete equation \eqref{B1} are of 
interest to several applications, for instance, in the study of
polymerization processes.

Coagulation equations with two-components, i.e., $d=2$, have been introduced
in \cite{L}. In that paper, the solutions to the coagulation equations with
kernels of the form $K(x_{1},x_{2};y_{1},y_{2})=K_{1}(x_{1}+x_{2};y_{1}%
+y_{2})$ have been written in terms of the solutions of the one-component
coagulation equation with kernel $K_{1}$. In particular, the solutions with
the constant kernel $K(x,y)=1$ are computed explicitly.

Multicomponent coagulation equations have been extensively used to analyse the
evolution of chemical properties of aerosol particles in atmospheric science
(cf.\ \cite{Olenius, Vehkam}). Additional details about the
physics of these systems can be found in \cite{FLNV}, \cite{FLNV2},
\cite{Fried}.

This paper is centred around a phenomenon that is specific to multicomponent
coagulation equations and that we termed in \cite{FLNV2} as
\emph{localization.} 
We say that the solutions to  equations (\ref{B1}), (\ref{B4}) localize if for large particle sizes they  tend to concentrate
along a  line, more precisely, along a ray starting from the origin $O$, in the space of compositions ($\mathbb{N}_{0}^{d}\setminus
\{O\}$ or ${\mathbb{R}}_{\ast}^{d}$, respectively).  Localization
can be observed, using generating functions, in the solutions of (\ref{B1})  for some particular kernels for which the solutions of the time dependent problem
(\ref{B1}) can be explicitly computed (cf.\ \cite{KBN} and also \cite{chp} for further details). 
We also remark that  the solutions of the continuous equation (\ref{B4}) can be computed in the case of the constant and the additive kernels using multicomponent Laplace transform  \cite{FDGG, L}.
In both cases the solutions  concentrate for long times along a
ray of the space of cluster concentrations $\mathbb{N}_{0}^{d}\setminus
\{O\}$, as discussed above. The orientation of
this ray is uniquely prescribed in terms of the initial distribution of
cluster compositions. 
More precisely, the conservation laws imply that the 
monomer composition of clusters, relative to the total number of monomers, along the ray follows the relative composition of monomers in the initial data, as will be explained below.

Interestingly, the localization direction is not encoded in any
property of the coagulation kernel. Indeed, we show in this paper that localization takes place for a large class
of coagulation kernels for which there is not any strongly preferred direction in the space of cluster compositions. Due to this we can think of localization as an emergent property.

The key relevance of  localization results is that they allow to reduce the long-time 
dynamics of multicomponent coagulation equations  to the
dynamics of just one-component coagulation systems. An example of how this
idea can be applied to specific cases is given later in Theorem
\ref{TheorConstLines}.

In this paper we  prove the localization of the solutions of the time
dependent problems (\ref{B1}), (\ref{B4}) for general classes of coagulation
kernels $K_{\alpha,\beta},\ K\left(  x,y\right)  .$ As mentioned above, the direction of localization is determined by the initial distribution of clusters. This could
be expected, because formally the following mass conservation properties hold
for (\ref{B1}) and (\ref{B4}) respectively:%
\begin{equation}
\partial_{t}\left(  \sum_{\alpha\in\mathbb{N}_{0}^{d}\setminus\{O\}}\alpha
n_{\alpha}\right)  =0\ \ ,\ \ \partial_{t}\left(  \int_{{\mathbb{R}}_{\ast
}^{d}}xf\left(  x,t\right)  dx\right)  =0  .\label{MassCons}%
\end{equation}
Notice that the identities in (\ref{MassCons}) are vector identities. This is
due to the conservation of the different types of monomers for the solutions
of (\ref{B1}) and (\ref{B4}). We will restrict our attention  to
coagulation kernels for which gelation does not take place. Therefore, the
identities (\ref{MassCons}) will be satisfied for the solutions considered in
this paper.  

Assuming the conservation laws and localization, we can also deduce the direction of localization.  For this, let us consider the standard $1$-norm 
$|\cdot|$ on ${\mathbb{R}}^{d}$, using which
$\left\vert x\right\vert =\sum_{j=1}^{d} x_{j}$ for $x\in \mathbb{R}^{d}_\ast$. Now, if $f$ at time $t$ is concentrated along a ray in some fixed direction $\theta_0$ with $|\theta_0|=1$, then to each $x$ in this region of concentration
$x\approx |x|\theta_0$.  But then the conserved total mass vector satifies
$m(t):=\int_{{\mathbb{R}}_{\ast
}^{d}}xf\left(  x,t\right)  dx \approx \theta_0
\int_{{\mathbb{R}}_{\ast
}^{d}} |x| f\left(  x,t\right)  dx$
which implies that $(\theta_0)_i \approx m_i(0)/|m(0)|$. Hence, for completely localized solutions the direction $\theta_0$ is already determined by the initial data.

The phenomenon of localization takes place also for generalizations of 
equations (\ref{B1}), (\ref{B4}) for which an additional source term is
included. In \cite{FLNV2} we provide a detailed study of the localization properties of
stationary solutions of (\ref{B1}), (\ref{B4}) with an additional source of
clusters on the right-hand side. More precisely, the equations studied in
\cite{FLNV2} are:%
\begin{align}
\kappa_d[n_\alpha] + s_{\alpha}  &
=0\,, \label{BB6}\\
\K_d[f](x)
+\eta\left(  x\right)   &  =0\,, \label{B6}%
\end{align}
where $s_{\alpha}\geq0,\ \eta\left(  x\right)  \geq0$ satisfy suitable
integrability conditions. Solutions to  equations \eqref{BB6}-\eqref{B6}
have been proven to exist in \cite{FLNV4} under general conditions on the
kernels $K_{\alpha,\beta},\ K\left(  x,y\right)  $ satisfying the conditions
\eqref{eq:condK_sym1}--\eqref{KernelRepres}. It has been proved in \cite{FLNV2}
that those solutions are concentrated along rays in the spaces $\mathbb{N}%
_{0}^{d}\setminus\{O\}$ and$\ {\mathbb{R}}_{\ast}^{d}$ for large
values of $\left\vert \alpha\right\vert $ and $\left\vert x\right\vert $,
respectively.
The direction of the localization line then depends on the first moments of the
source terms $s_{\alpha}$ and $\eta\left(  x\right)  .$

In this paper we  prove that localization takes place for the time dependent solutions of (\ref{B1}) and (\ref{B4}) for a large class of coagulation kernels. In the case of equation (\ref{B4}) we assume that
\begin{equation}
K\in C(({\mathbb{R}}_{*}^{d})^{2}),\ K(x,y)=K(y,x),\ K(x,y)\geq0.
\label{eq:condK_cont}%
\end{equation}
We require continuity of the kernels $K$ in order to obtain meaningful
formulas for measure-valued solutions $f.$ 
In addition, we will assume that
\begin{equation}
K(x,y)\geq c_{1}(|x|+|y|)^{\gamma}\Phi_{p}\left(  \frac{|x|}{|x|+|y|}\right)
\ \ ,\ \ x,y\in{\mathbb{R}}_{\ast}^{d} \label{eq:condK_sym1}%
\end{equation}%
\begin{equation}
K(x,y)\leq c_{2}(|x|+|y|)^{\gamma}\Phi_{p}\left(  \frac{|x|}{|x|+|y|}\right)
\ \ ,\ \ x,y\in{\mathbb{R}}_{\ast}^{d} \label{eq:condK_sym2}%
\end{equation}
with $\gamma\in{\mathbb{R}}$, and for some 
$p\in{\mathbb{R}}$ and 
$\Phi_{p}\in C(0,\infty)$ such
that
\begin{equation}
\Phi_{p}(s)=\frac{1}{s^{p}(1-s)^{p}}%
,\quad0<s<1,\text{ with }\gamma+2p\geq 0\,, \label{eq:Phi_kernel}%
\end{equation}
and some constants $0<c_{1}\leq c_{2}<\infty$. 
Note that then $\Phi_{p}(s)=\Phi_{p}(1-s)$
and thus the bounds are symmetric functions, due to
$\frac{|y|}{|x|+|y|}=1-\frac{|x|}{|x|+|y|}$.

In the case of the equation
(\ref{B1}), we will assume that the kernel $K_{\alpha,\beta}$ can be written
as%
\begin{equation}
K_{\alpha,\beta}=K(\alpha,\beta)\ \ \text{for }\alpha,\beta\in\mathbb{N}%
_{0}^{d}\setminus\{O\} \label{KernelRepres}%
\end{equation}
for some $K(x,y)$ satisfying the earlier requirements. Notice that if $K_{\alpha,\beta}$ 
is a kernel for which there are 
$c_1,c_2$ such that 
$c_{1}(|\alpha|+|\beta|)^{\gamma}%
\Phi_{p}(\frac{|\alpha|}{|\alpha|+|\beta|})\leq K_{\alpha,\beta}\leq
c_{2}(|\alpha|+|\beta|)^{\gamma}\Phi_{p}(\frac{|\alpha|}{|\alpha|+|\beta|})$
for all $\alpha,\beta\in\mathbb{N}_{0}^{d}\setminus\{O\}$, we may find a function $K(x,y)$ of the required type so that  (\ref{KernelRepres}) holds (changing, if needed, the values of the constants $c_{1},c_{2}$). We
remark that the estimates \eqref{eq:condK_sym1}--\eqref{eq:Phi_kernel} are
invariant under the permutation of the components $x_{1},x_{2},...,x_{d}$. In
particular, the kernels satisfying \eqref{eq:condK_sym1}--\eqref{eq:Phi_kernel}
cannot have different power law behaviour along any two different variables,
say $x_{j},\ x_{k}$ with $j\neq k$.

In order to avoid gelation, we will assume also the following conditions for
the parameters $\gamma$ and $p$ in (\ref{eq:condK_sym1}), (\ref{eq:condK_sym2}%
), (\ref{eq:Phi_kernel}):
\begin{equation}
 \gamma+p<1,\ \ \ \gamma<1 .\label{KernIneq}
\end{equation}
We would like to point out that the class of kernels considered here strictly
contains the class commonly found in the literature, namely, those satisfying the power law
bounds
\begin{equation}
c_1(|x|^{\gamma+\lambda}|y|^{-\lambda}+|y|^{\gamma+\lambda}|x|^{-\lambda})\leq
K(x,y)\leq c_2(|x|^{\gamma+\lambda}|y|^{-\lambda}+|y|^{\gamma+\lambda}%
|x|^{-\lambda}) \,, \label{eq:kernel_lambda}%
\end{equation}
\[
-\lambda<1,\quad\gamma+\lambda<1\ \ ,\ \ \gamma<1\, .
\]
Then we may choose $p=\max\{\lambda,-\gamma-\lambda\}$, for which $\gamma+p\geq-p$ due to
$\gamma+2p=|\gamma+2\lambda|$. We recall that the coagulation kernels depend
on the specific mechanism which is responsible for the aggregation of the
clusters at the microscopic level. In particular, the class of kernels
\eqref{eq:kernel_lambda} contains the physically relevant kernels that are often
used in aerosol science, such as the free molecular (ballistic) kernel and the Brownian kernel  (cf.\ \cite{Fried, Vehkam} as well as \cite{ FLNV2} for a more detailed discussion). 

Notice that the kernels $K$ satisfying (\ref{eq:condK_sym1}),
(\ref{eq:condK_sym2}), (\ref{eq:Phi_kernel}) are bounded from above and below
by homogeneous functions, but they are not necessarily themselves homogeneous. In some of
the results presented later we will need to assume homogeneity, i.e., then we additionally require that
\begin{equation}
K(rx,ry)=r^{\gamma}K(x,y)\ \ \text{,\ \ }r>0\ \ ,\ \ x,y\in{\mathbb{R}}_{\ast
}^{d}\,. \label{eq:condK_homogeneity}%
\end{equation}

In the one component case ($d=1$),
there are already many earlier results about solutions  
to  equations (\ref{B1}) and (\ref{B4}) available.
For example,
assuming that (\ref{eq:condK_sym1}), (\ref{eq:condK_sym2}),
(\ref{eq:Phi_kernel}), (\ref{KernIneq}) are satisfied
and that the
initial mass ($\sum_{\alpha\in\mathbb{N}_{0}^{d}\setminus\{O\}}\alpha
n_{\alpha}$ or $\int_{{\mathbb{R}}_{\ast}^{d}}xf\left(  x,t\right)  dx$,
respectively) is finite, then the mass becomes concentrated in
the region of cluster sizes of order $\alpha\approx t^{\frac{1}{1-\gamma}}$ or
$x\approx t^{\frac{1}{1-\gamma}}$, respectively, as $t\rightarrow\infty$.
Moreover, if $d=1$, $0\leq\gamma<1$, $0\leq\gamma+p<1$, and
(\ref{eq:condK_homogeneity}) holds, it is well known that self-similar
solutions of (\ref{B4}) with the form $f\left(  x,t\right)  =\frac{1}%
{t^{\frac{2}{1-\gamma}}}F\left(  \frac{x}{t^{\frac{1}{1-\gamma}}}\right)  $
exist (cf.\ \cite{EM05,EM06,FL04}). These self-similar
solutions are expected to represent the long time asymptotics of the solutions of
(\ref{B1}) and (\ref{B4}) in great generality, although this has been rigorously proven only for 
 particular kernels from the class defined by (\ref{eq:condK_sym1}), (\ref{eq:condK_sym2}), (\ref{eq:Phi_kernel}), (\ref{KernIneq}), specifically, only if the kernel $K$ is
constant \cite{MP04} or a perturbation of a constant \cite{CT21, Seb2}. 
It is also possible to obtain  representation formulas using Laplace transforms with the additive kernel $K(x,y) = x+y$ and with the multiplicative kernel $K(x,y) = x  y$ (cf.\ \cite{MP04}, \cite{MP06}). In the multicomponent case ($d \geq 1$),  representation formulas 
for the solutions of the initial value problem associated to (\ref{B1}) (or
(\ref{B4})) 
can also be obtained for the constant kernel $K(x,y)=1$, the additive kernel $K(x,y)=|x|+|y|$ and the
product kernel $K(x,y)=|x||y|$ (cf.\ \cite{FDGG, FDGG2, KBN, L}) using multicomponent Laplace transform methods.

Analogous estimates, which show that the mass of the clusters is concentrated
in the self-similar region (i.e. $\left\vert \alpha\right\vert \approx
t^{\frac{1}{1-\gamma}}$or $\left\vert x\right\vert \approx t^{\frac
{1}{1-\gamma}}$) for large times $t$, can be derived in the case of
multicomponent coagulation systems, adapting in a suitable manner the methods
used to prove these results in the case $d=1$. We remark that the total mass
of monomers, defined as $\sum_{\alpha\in\mathbb{N}_{0}^{d}\setminus
\{O\}}\left\vert \alpha\right\vert n_{\alpha}$ or $\int_{{\mathbb{R}}_{\ast
}^{d}}\left\vert x\right\vert f\left(  x,t\right)  dx$, respectively for
(\ref{B1}) and (\ref{B4}),  is conserved (cf.\ (\ref{MassCons})). In
addition to such estimates, we also prove that the mass is concentrated along a
particular ray of the cluster space 
as
$t\rightarrow\infty,$ i.e.,  that localization in the sense defined above
takes place. The localization results for time dependent problems are the
main novelty of the present paper. Notice that, since the solutions of (\ref{B1}) can
be interpreted as particular solutions of (\ref{B4}) with the form (\ref{B5}),
the localization results for (\ref{B1}) will follow from the
corresponding results for (\ref{B4}).

As a final remark, we note that, except for particular kernels such as the constant kernel \cite{MP06} or kernels that are 'close' to constant \cite{ NTV16, Seb2}, in general there  are no uniqueness results for self-similar solutions available in the literature.
Nevertheless, our result imply that all self-similar solutions localize. Notably, for a special class of kernels that are constant only along rays, it is possible to obtain uniqueness and stability results by employing ealier results for the one-component equation with constant kernel (cf.\ Theorem
\ref{TheorConstLines}).

\subsection{Notations}

We collect here for the reader's convenience the main notations and definitions which will be repeatedly used throughout the paper,
some of these having been already introduced above.

We denote ${\mathbb{R}}_{+}=[0,\infty)$, to be distinguished from the already defined sets ${\mathbb{R}}_{\ast}=(0,\infty)$ and ${\mathbb{R}}_{\ast}^{d}%
=[0,\infty)^{d}\backslash\{O\}$.
We use $|\cdot|$ and $\Vert\cdot\Vert$ to denote the following norms on
${\mathbb{R}}^{d}$, 
\[
\left\vert x\right\vert =\sum_{j=1}^{d}\left\vert x_{j}\right\vert\,, \quad
\left\Vert x\right\Vert =\sqrt{\sum_{j=1}^{d}\left(  x_{j}\right)
^{2}}\,, \qquad\text{for } x\in{\mathbb{R}}_{\ast}^{d},\ x=\left(  x_{1}%
,x_{2},...,x_{d}\right).
\]

We denote by $C_{c}\left(  {{\mathbb{R}}_{\ast}^{d}}\right)  $ the set of
compactly supported continuous functions in ${{\mathbb{R}}_{\ast}^{d}}$,
and by 
$C_{c}^{k}\left(  {{\mathbb{R}}_{\ast}^{d}}\right)$, for $k=1,2,\ldots$, 
the set of such compactly supported functions with
$k$ continuous derivatives.
We use the notations $\mathcal{M}_{+}\left(  {\mathbb{R}}_{+}\right)  $
and $\mathcal{M}_{+}\left(  {\mathbb{R}}_{\ast}^{d}\right)  $ to denote the
 spaces of non-negative   Radon measures on ${\mathbb{R}}_{+}$ and ${\mathbb{R}%
}_{\ast}^{d}$, respectively. We will use indistinctly the notation $f\left(
dx\right)  $, $f\left(  x\right)  dx$, or $f$ to denote a measure
$f\in\mathcal{M}_{+}\left(  {\mathbb{R}}_{+}\right)  $ or in $\mathcal{M}%
_{+}\left(  {\mathbb{R}}_{\ast}^{d}\right)  .$ The former notation will be
preferred when the measure is integrated against a test function. We stress
that we will use the notation $f\left(  x\right)  dx$ to denote a measure on
$\mathbb{R}_{\ast}^{d}$ even if this measure is not absolutely continuous with
respect to the Lebesgue measure of $\mathbb{R}_{\ast}^{d}.$

We denote by $\Delta^{d-1}$ the simplex
\begin{equation}
\Delta^{d-1}=\left\{  \theta\in\mathbb{R}_{\ast}^{d}:\left\vert \theta
\right\vert =1\right\}  \label{eq:simplex}%
\end{equation}
and as $\mathcal{M}_{+}\left(  \Delta^{d-1}\right) , \mathcal{M}_{+}\left(
{\mathbb{R}}_{*} \times\Delta^{d-1}\right)  $ the spaces of non-negative Radon
measures on $\Delta^{d-1}$ and on ${\mathbb{R}}_{*} \times\Delta^{d-1}$,
respectively. We will denote by 
$\delta_{\theta_{0}}$ or
$\delta\left(  \cdot-\theta_{0}\right)  $ the
Dirac measure supported at $\theta_{0}\in\Delta^{d-1}.$

We will denote by $C$ a generic constant which can depend on $d$ and on the properties
of the kernels (specifically, $\gamma$, $p$, as well as $c_{1}$ and $c_{2}$ in
(\ref{eq:condK_sym1}), (\ref{eq:condK_sym2})) but which is independent of the solution
under consideration. The value of $C$ may also change from line to line.

\subsection{Main results}

We now state the main results proved in this paper. The precise definitions will be given later in Section \ref{DefSect}.  We begin with our main localization
result.

\begin{theorem}
\label{mainRoughImproved} Let $f_0 \in \mathcal{M}_+(\R_*^d)$ satisfying $\int_{{\mathbb{R}}_{\ast}^{d}
}(\left\vert x\right\vert + \vert x\vert ^{1+r})f_{0}\left(  dx\right)  <\infty$ for some
$r>0$. Define $m(0):= \int_{{\mathbb{R}}_{\ast}^{d}}x f_0\left( d x\right)\in \R^d$, denote $m_0=|m(0)|$, and suppose that $m_0>0$. 
Let the coagulation kernel $K$ satisfy the assumptions (\ref{eq:condK_cont}), (\ref{eq:condK_sym1}),
(\ref{eq:condK_sym2}), (\ref{eq:Phi_kernel})  with $0\leq\gamma<1$, 
and   $0\leq\gamma+p<1.$ 
Then there exists a weak solution $f\in C\left(  \left[  0,\infty\right)
;\mathcal{M}_{+}\left(  {\mathbb{R}}_{\ast}^{d}\right)  \right)  $ to (\ref{B4}), (\ref{coagKern}) such that $f(\cdot,0) = f_0$  with the following properties.  This solution is mass-conserving:
$\int_{{\mathbb{R}}_{\ast}^{d}}
x f\left(  x,t\right)  dx=m(0)$ for all $t\geq 0$, and in addition it satisfies 
\begin{equation}
\int_{{\mathbb{R}}_{\ast}^{d}}\left\vert x\right\vert ^{k}f\left(  x,t\right)
dx\leq C_{0}t^{\frac{k-1}{1-\gamma}}\,, \quad
t\geq1\,, \label{eq:moment_estimate}%
\end{equation}
 for some $k>1$ and $C_{0}>0.$
Moreover, there exists a function $\delta\left(  \cdot\right)  \in C\left(
\left[  0,\infty\right)  \right)  $ such that $\delta\left(  t\right)  >0$ for
$t\in\left[  1,\infty\right)  $ and $\lim_{t\rightarrow\infty}\delta\left(
t\right)  =0$ and for which
\begin{equation}
\lim_{t\rightarrow\infty}\left\vert \int_{\left\{  \delta\left(  t\right)
t^{\frac{1}{1-\gamma}}\leq\left\vert x\right\vert \leq\left(  \delta\left(
t\right)  \right)  ^{-1}t^{\frac{1}{1-\gamma}}\right\}  \cap\left\{
\left\vert \frac{x}{\left\vert x\right\vert }-\theta_{0}\right\vert \leq
\delta\left(  t\right)  \right\}  }\left\vert x\right\vert f\left(
x,t\right)  dx-m_{0}\right\vert =0 \label{localRestricted}%
\end{equation}
where%
\begin{equation}
\theta_{0}=\frac{\int_{{\mathbb{R}}_{\ast}^{d}}x f_0\left(  x\right)  dx}
{m_{0}}\in \Delta^{d-1}\ \ ,\ \ \left\vert \theta_{0}\right\vert =1.
\label{DefThetaZero}%
\end{equation}
\end{theorem}

\begin{remark}
The crucial information about the function $\delta\left(  t\right)  $ is that
it converges to zero. Therefore, (\ref{localRestricted}) implies that the mass
is localized along a particular direction in distances $\left\vert
x\right\vert $ of order $t^{\frac{1}{1-\gamma}}$ for all the large times $t$
as $t\rightarrow\infty$, i.e., localization of the measure $\left\vert
x\right\vert f\left(  x,t\right)  dx$ takes place as $t\rightarrow\infty.$
Notice that the vector $\theta_{0}$ defined in (\ref{DefThetaZero}) only depends on the conserved quantities, just as was discussed in the Introduction.
\end{remark}

Notice that Theorem \ref{mainRoughImproved} yields localization for a particular weak solution
of the initial value problem (\ref{B4}), (\ref{coagKern}) with initial value $f(0,.)=f_0$. The reason why the localization result is not stated for  every weak solution  is due to the lack of  a uniqueness theory. Indeed, the arguments used in the proof of  Theorem \ref{mainRoughImproved} rely on the results of \cite{EM06} that only ensure existence of a weak solution to (\ref{B4}), (\ref{coagKern}) satisfying \eqref{eq:moment_estimate} with initial value $f_0$, but no uniqueness is proved in \cite{EM06}. 
A theory of uniqueness of weak solutions  combined with Theorem \ref{mainRoughImproved} would then imply localization for all  weak solutions of (\ref{B4}), (\ref{coagKern}).
 The derivation of such  results for weak solutions is not the goal of this paper. 
Uniqueness results in the one-component case $d=1$, for some kernels satisfying the upper bound \eqref{eq:condK_sym2} with $\gamma \leq 1 $ and $p=0$ as well as  additional regularity conditions have been obtained in \cite{FL06}.

 The condition \ \eqref{eq:moment_estimate}    ensures that most of the mass of the solution remains in the self-similar region. We expect the estimate \eqref{eq:moment_estimate} to hold for all  weak solutions to (\ref{B4}), (\ref{coagKern}) which decay sufficiently fast for large $|x|$ and for the  physically relevant kernels with homogeneity smaller than $1$.
It turns out that  it is possible to obtain a slightly weaker localization result  for all  solutions to (\ref{B4}), (\ref{coagKern}) satisfying the moment estimate \eqref{eq:moment_estimate} for a more general class of kernels than the one considered in Theorem \ref{mainRoughImproved}.
 More precisely we have the following result.

\begin{theorem}
\label{mainRough}Suppose that $f\in C\left(  \left[  0,\infty\right)
;\mathcal{M}_{+}\left(  {\mathbb{R}}_{\ast}^{d}\right)  \right)  $ is a weak
solution of (\ref{B4}) such that $0<\int_{{\mathbb{R}}_{\ast}^{d}}\left\vert
x\right\vert f\left(  x,t\right)  dx=m_{0}<\infty$ and such that the
assumptions \eqref{eq:condK_cont}, (\ref{eq:condK_sym1}), (\ref{eq:condK_sym2}), (\ref{eq:Phi_kernel}) hold with {$\gamma{,\ p}$ satisfying} (\ref{KernIneq}). Assume that 
there are $a>1$ and $C_{0}>0$ such that $f$
satisfies (\ref{eq:moment_estimate}) 
for all $k\in
\lbrack1/a,a]$. Then, there exists a function
$\delta\left(  \cdot\right)  \in C\left(  \left[  0,\infty\right)  \right)  $
such that $\delta\left(  t\right)  >0$ for $t\in\left[  1,\infty\right)  $ and
$\lim_{t\rightarrow\infty}\delta\left(  t\right)  =0$ as well as a Borel set
$I\subset\left[  0,\infty\right)  $ with the property that $\lim
_{T\rightarrow\infty}\frac{\left\vert I\cap\left[  T,2T\right]  \right\vert
}{T}=0$ such that
\begin{equation}
\lim_{T\rightarrow\infty}\left(  \sup_{t\in\left[  T,2T\right]  \backslash
I}\left\vert \int_{\left\{  \delta\left(  t\right)  t^{\frac{1}{1-\gamma}}%
\leq\left\vert x\right\vert \leq\left(  \delta\left(  t\right)  \right)
^{-1}t^{\frac{1}{1-\gamma}}\right\}  \cap\left\{  \left\vert \frac
{x}{\left\vert x\right\vert }-\theta_{0}\right\vert \leq\delta\left(
t\right)  \right\}  }\left\vert x\right\vert f\left(  x,t\right)
dx-m_{0}\right\vert \right)  =0 \label{locEstimate}%
\end{equation}
where $\theta_{0}$ is as in (\ref{DefThetaZero}).
\end{theorem}

Notice that the main difference between Theorems  \ref{mainRoughImproved} and \ref{mainRough}
 is that in the first case we assume a more restricted
set of parameters $\gamma$ and $p.$ On the other hand, we obtain stronger
localization results in the case of Theorem \ref{mainRoughImproved}. In  (\ref{locEstimate}) we allow for the existence of a set of times
$I\subset\left(  1,\infty\right)  $ whose density converges to zero for large
values of $t$ and for which the localization property could fail. We do not
know if it is possible to have solutions of (\ref{B4}) for which localization
does not take place for a small set of large times. Most likely such a type
of behaviour does not take place for any solution of (\ref{B4}). However, only
the estimate (\ref{locEstimate}) can be obtained from the assumptions on
the solutions to (\ref{B4}) considered in Theorem \ref{mainRough}. In fact, estimates ensuring that the mapping $t \mapsto (1+t)^{\frac{d+1}{1-\gamma}}f((t+1)^{\frac{1}{1-\gamma}}\cdot,t)$ is uniformly continuous in the weak-$\ast$ topology
would yield the
stronger localization result (\ref{localRestricted}), 
but this would require analysis going much beyond the currently available well posedness results.

We will discuss in Sections \ref{GlobEx} and \ref{MomEst} sufficient
conditions for the moment estimates \eqref{eq:moment_estimate} to be satisfied
for the range of parameters $0\leq\gamma<1,\ 0\leq\gamma+p<1$, which  in particular  is contained in the range defined by
\eqref{KernIneq}. In the case of moments $k>1$ we need to assume suitable
conditions on the initial data $f_{0}.$

Since the solutions to the discrete coagulation equation (\ref{B1}) are
particular solutions of (\ref{B4}) having the form (\ref{B5}) the localization
result in Theorem \ref{mainRoughImproved} holds for the solutions of
(\ref{B1}). Given that this result has an independent interest, we formulate
it here separately. (It would be possible to formulate also a discrete version
of Theorem \ref{mainRough}).

\begin{theorem}
\label{thm:concentDisc} Suppose that $K_{\left(  \cdot,\cdot\right)  }:\left(
\mathbb{N}_{0}^{d}\setminus\{O\}\right)  ^{2}\rightarrow\mathbb{R}_+$ is a
mapping that can be written in the form (\ref{KernelRepres}) for some function
$K:({\mathbb{R}}_{\ast}^{d})^{2}\rightarrow{\mathbb{R}}_{+}$ satisfying
(\ref{eq:condK_cont}) as well as the bounds
\eqref{eq:condK_sym1}--\eqref{eq:Phi_kernel} for some $\gamma\in\lbrack0,1)$
and $p\in\mathbb{R}$ such that $0\leq\gamma+p<1.$ 
Let $\{n_{\alpha,0}\}_{{\alpha\in\mathbb{N}_{0}^{d}\setminus\{O\}}}$ satisfy $0<\sum_{\alpha\in\mathbb{N}_{0}^{d}\setminus\{O\}}|\alpha|^{1+r}n_{\alpha,0}  <\infty$ for some $r>0$.
Then there is a solution  $\left\{  n_{\alpha
}\left(  \cdot\right)  \right\}  _{\alpha\in\mathbb{N}_{0}^{d}\setminus
\{O\}}$  of \eqref{B4} such that $n_\alpha(0) = n_{\alpha,0}$ and, for each $t>0$,
\[
\sum_{\alpha\in\mathbb{N}_{0}^{d}\setminus\{O\}}\left\vert \alpha\right\vert
n_{\alpha}(t) = \sum_{\alpha\in\mathbb{N}_{0}^{d}\setminus\{O\}}\left\vert \alpha\right\vert
n_{\alpha,0}:=m_{0} \in (0,\infty).
\]
Moreover, it satisfies the following localization property. There exists a positive function $\delta\in C(0,\infty)$ such that
$\lim_{t\rightarrow\infty}\delta(t)=0$ and with the property that 
\begin{equation}
\lim_{t\rightarrow\infty}\left\vert \sum_{\left\{  \delta\left(  t\right)
t^{\frac{1}{1-\gamma}}\leq\left\vert \alpha\right\vert \leq\left(  \delta\left(
t\right)  \right)  ^{-1}t^{\frac{1}{1-\gamma}}\right\}  \cap\left\{
\left\vert \frac{\alpha}{\left\vert \alpha\right\vert }-\theta_{0}\right\vert \leq
\delta\left(  t\right)  \right\}  }\left\vert \alpha\right\vert n_{\alpha
}\left(  t\right)  -m_{0}\right\vert =0 \label{LocalDisc}%
\end{equation}
where $\theta_{0}\in \Delta^{d-1}$ is defined by means of
\begin{equation}
\theta_{0}:=\frac{\sum_{\alpha\in\mathbb{N}_{0}^{d}\setminus\{O\}}\alpha
n_{\alpha,0}}{m_{0}}. \label{LocDirectDisc}%
\end{equation}
\end{theorem}

We will also study localization properties for the self-similar solutions of
(\ref{B4}) with $d>1.$ The mass conserving self-similar solutions are
solutions of (\ref{B4}) with the form:
\begin{equation}
f\left(  x,t\right)  =(\varepsilon_t) ^ {1+d} F\left( x \varepsilon_t \right)
\ \ ,\ \ \xi=x\varepsilon_t\ \ , \ \ \varepsilon_t = (t+1)^{-\frac{1}{1-\gamma}}.
\label{MultSelfSimForm}%
\end{equation}
The existence of solutions of (\ref{B4}) with the form (\ref{MultSelfSimForm})
under the assumptions (\ref{eq:condK_sym1}), (\ref{eq:condK_sym2}),
(\ref{eq:Phi_kernel}), (\ref{KernIneq}), (\ref{eq:condK_homogeneity}) and for
$0\leq\gamma<1$, $0\leq\gamma+p<1$ has been proved in the case $d=1$ in
\cite{EM05,EM06,FL04}. Using these results it is possible to
prove the existence of self-similar solutions in the multicomponent case $d>1$
under analogous assumptions on the collision kernels and having the particular
form
\begin{equation}
F\left(  \xi\right)  = \frac{\sqrt{d}}{|\xi|^{d-1}} F_{0}\left(  \left\vert \xi\right\vert \right)
\delta\left(  \frac{\xi}{\left\vert \xi\right\vert }-\theta_{0}\right)
\label{SelfSimLoc}
\end{equation}
where $\theta_{0}\in\Delta^{d-1}$, $\delta \in \mathcal{M}_+(\Delta^{d-1})$ is supported at $\theta_0$ and $F_{0}$ is a self-similar profile for a suitable
one-dimensional coagulation equation. The existence of self-similar profiles
with the form (\ref{SelfSimLoc}) will be seen in Section \ref{GlobEx}.

It turns out that all the solutions of (\ref{B4}) with the form
(\ref{MultSelfSimForm}) and satisfying suitable integrability conditions for both 
small and large $\left\vert \xi\right\vert $, have
the form (\ref{SelfSimLoc}). This result can be interpreted as a localization
result analogous to the Theorems \ref{mainRoughImproved}, \ref{mainRough} for
solutions of (\ref{B4}). The precise localization result for self-similar
solutions that we will prove in this paper is the following.

\begin{theorem}
\label{TheorSelfSimLoc}Suppose that the assumptions (\ref{eq:condK_sym1}),
(\ref{eq:condK_sym2}), (\ref{eq:Phi_kernel}), (\ref{KernIneq}) are satisfied.
Suppose that $F\in\mathcal{M}_{+}\left(  {\mathbb{R}}_{\ast}^{d}\right)  $ is
a self-similar profile with finite mass for (\ref{B4}) in the sense of
Definition \ref{def:self_sim_sol}. Then, there exists $\theta_{0}%
\in \Delta^{d-1}$  such
that $F$ has the form (\ref{SelfSimLoc}) where $F_{0}$ is a self-similar
profile associated to the one component coagulation equation (i.e. $d=1$) and
coagulation kernel $K_{\theta_0}\left(  s,r\right)  =K\left(  s\theta_{0}%
,r\theta_{0}\right)  ,\ s,r\in{\mathbb{R}}_{\ast}.$
\end{theorem}

An interesting consequence of the localization results contained in Theorem
\ref{mainRoughImproved} is that they allow to characterize the long time
asymptotics for a class of coagulation kernels for which it does not seem
feasible to obtain an explicit representation formula for the solutions. We
recall that in the case of one component systems a complete characterization
of the long time asymptotics for arbitrary initial data has been obtained only
for coagulation kernels with homogeneity smaller than one for which it is possible to obtain representation
formulas of the solutions using Laplace transform methods (cf.\ \cite{MP04,
MP06}), or for kernels $K$ which are close to the constant kernel (cf.
\cite{Seb2}).

We will combine the localization results obtained in this paper (cf.\ Theorem
\ref{mainRoughImproved}) with the characterization of the long time
asymptotics obtained in \cite{MP04, MP06} to characterize the long time
behaviour of the solutions of coagulation equations with kernels that are
constant along each ray that passes through the origin. This is due
to the fact that the localization of the solutions along a ray allows to
approximate the behaviour of the solutions by a one-component coagulation
equation with a constant kernel. It does not seem feasible to derive an
explicit formula using Laplace transform methods for the class of kernels with
the form (\ref{KLines}) below, except for some very particular choices of the
function $Q.$ We have the following result.

\begin{theorem}
\label{TheorConstLines}Suppose that the kernel $K$ satisfies \eqref{eq:condK_cont} and has the form
\begin{equation}
K(r\theta,s\theta)=Q(\theta) \label{KLines}%
\end{equation}
for any $r,s>0$ and for any $\theta=\left(  \theta_{1},\theta_{2}%
,...,\theta_{d}\right)  \in\Delta^{d-1}$. Here $Q$ is a continuous function defined on $\Delta^{d-1}  $ and
$0<c_{1}\leq Q\left(  \theta\right)  \leq c_{2}.$ Let $f_{0}\in L^{1}\left(
\mathbb{R}_{\ast}^{d}\right)  $ be a nonnegative function satisfying
\begin{equation}
m_0:= \int_{\mathbb{R}_{\ast}^{d}}|x|f_{0}\left(  x\right)  dx >0 \label{F0Mass}%
\end{equation}
and also%
\begin{equation}
\int_{\mathbb{R}_{\ast}^{d}}\left\vert x\right\vert^{a}f_{0}\left(  x\right)  dx<\infty\label{MomHomOne}%
\end{equation}
for some $a>1.$ Then there exists a function $f\in C\left(  \left[  0,\infty
\right)  ;L^{1}\left(  \mathbb{R}_{\ast}^{d}\right)  \right)  \cap 
C^{1}\left(  \left(  0,\infty\right)  ;L^{1}\left(  \mathbb{R}_{\ast}%
^{d}\right)  \right)  $ that solves (\ref{B4}) in the classical sense, satisfying $f\left(  x,0\right)  =f_{0}\left(  x\right) $  and 
$$\int_{\R^d_*} xf(x,t)dx = m:= \int_{\R^d_*} xf_0(x)dx,\ t>0.$$ Moreover, we have%
\[
\lim_{t\rightarrow\infty}t^{2}f\left(  t\xi,t\right)  =F_{0}\left(  \left\vert
\xi\right\vert ;\theta_{0}\right)  \delta\left(  \frac{\xi}{\left\vert
\xi\right\vert }-\theta_{0}\right)
\]
where the convergence takes place in the weak$-\ast$ topology of
$\mathcal{M}_{+}\left(  \mathbb{R}_{\ast}^{d}\right)  $ and  where
$$\theta_{0}:=\frac{m}{m_0 }\in \Delta^{d-1} \quad \text{ and }$$
\[
F_{0}\left(  \left\vert \xi\right\vert ;\theta_{0}\right)  :=\frac{4\sqrt{d}%
}{\left(  Q(\theta_{0})\right)  ^{2}m_0 }\frac
{1}{\left\vert \xi\right\vert ^{d-1}}\exp\left(  -\frac{2\left\vert
\xi\right\vert }{Q(\theta_{0})m_0 }\right)  .
\]

\end{theorem}

Notice that the mass vector $\int_{\mathbb{R}_{\ast}^{d}}f\left(  x,t\right)
xdx=m$ remains constant for arbitrary values of $t\geq0.$ Theorem
\ref{TheorConstLines} states that for each value of $m\in\mathbb{R}_{\ast}%
^{d}$ there exists a unique self-similar solution of the form
(\ref{MultSelfSimForm}), (\ref{SelfSimLoc}) with $\theta_{0}=\frac
{m}{m_0 }$ which is a global attractor for the solutions
of (\ref{B4}) satisfying $m=\int_{\mathbb{R}_{\ast}^{d}}f_{0}\left(  x\right)
xdx$ as well as the moment estimate (\ref{MomHomOne}).

It seems possible to extend Theorem \ref{TheorConstLines} to initial values
$f_{0}$ in some measure spaces. A technical problem which arises if we try to
replace the space $L^{1}\left(  \mathbb{R}_{\ast}^{d}\right)  $ by the space
$\mathcal{M}_{+}\left(  \mathbb{R}_{\ast}^{d}\right)  $ is that the kernels
$K$ with the form (\ref{KLines}) are not necessarily continuous at $x=y=0$ and
therefore it is not possible to define the products $K\left(  x,y\right)
f\left(  dx\right)  f\left(  dy\right)  .$ In order to avoid these
technicalities we prefer to use the space $L^{1}\left(  \mathbb{R}_{\ast}%
^{d}\right)  .$

\subsection{Plan of the paper}

The plan of this paper is the following. In Section 2 we introduce several
definitions and notation that will be used in the rest of the paper. In Section 3 we prove the localization results for
the time dependent solutions and for the self-similar solutions. Specifically,
we prove Theorems  \ref{mainRoughImproved}, \ref{mainRough}, and
\ref{TheorSelfSimLoc}. The proof of these results is based on the use of some
particular test functions that are reminiscent of those used in the proof of
the localization results for stationary solutions in \cite{FLNV2}. 
 A difference with the arguments in \cite{FLNV2} is that in the situation considered in this paper it can be proven that most of the mass of the clusters  concentrate in the self-similar region, 
$\left\vert x\right\vert \approx t^{\frac{1}{1-\gamma}}$. 
This provides a natural cutoff for the solutions which allows to show that the contribution of the regions $x \ll  t^{\frac{1}{1-\gamma}}$ and 
$ x \gg t^{\frac{1}{1-\gamma}} $ is negligible. 
On the contrary, in the stationary solutions treated in  \cite{FLNV2} there is no characteristic cluster size in which most of the mass of the solutions is concentrated.

Section \ref{GlobEx}
collects several well  posedness results and moment estimates for the solution
of the coagulation equation (\ref{B4}). These results are well known in the
case of one component coagulation systems and their proof can be readily
adapted to the multicomponent case. The results in Section \ref{GlobEx} show that the
assumptions made on the solutions of the coagulation equations in Theorems
\ref{mainRoughImproved} and \ref{mainRough} hold for a suitable set of
parameters $\gamma$ and $p$ and a large class of initial data. We prove also
in this Section that the measures with the form (\ref{SelfSimLoc}) yield
self-similar solutions of the multicomponent coagulation equation if we assume
that $F_{0}$ is a self-similar profile of a suitable coagulation equation.
Section \ref{MomEst} contains the proof of certain moment estimates which constitute some of the 
key assumptions on the solutions necessary in order to prove the localization
results. These estimates prove that the mass of the solutions remain within
the self-similar region, $\left\vert x\right\vert \approx t^{\frac{1}%
{1-\gamma}}$ as $t\rightarrow\infty$, for a large class of initial data.
Although they are well known in the case of one-component coagulation systems,
these estimates play a crucial role in the proof of our localization
results, and we have written in detail the way in which their proof can be adapted
to the multicomponent coagulation case. In Section 6 we study the long time
asymptotics to the solutions of the multicomponent coagulation equation with
kernels satisfying (\ref{KLines}). In particular, the proof of Theorem
\ref{TheorConstLines} is given in this Section.

\section{Definitions and auxiliary results\label{DefSect}}

In this Section we provide the definition of weak solutions and self-similar
profiles that will be used in the following. We also collect, without proof,
several results for the multicomponent coagulation equation (\ref{B4}) that
are well known for one-component coagulation systems and can be proved for
multicomponent coagulation systems by means of simple adaptations of the
methods used to derive them in the one-component case.

We now introduce the definitions of solutions to (\ref{B1}), (\ref{B4}). We
formulate the definition of solution in the continuous case (\ref{B4}) since
the discrete case (\ref{B1}) can be considered as a particular case of
solutions $f$ having the form (\ref{B5}).

\begin{definition}
\label{def:time-dep_sol}Let $K$ be as in (\ref{eq:condK_cont}) and satisfy the
upper bound \eqref{eq:condK_sym2}, \eqref{eq:Phi_kernel}. Suppose that
\thinspace$f_{0}\in\mathcal{M}_{+}({{\mathbb{R}}_{\ast}^{d}})$ satisfies%
\begin{equation}
\int_{{\mathbb{R}}_{\ast}^{d}}|x|f_{0}(dx)<\infty . \label{CondInVal}%
\end{equation}
A function $f\in C([0,\infty);\mathcal{M}_{+}({{\mathbb{R}}_{\ast}^{d}}))$
with $\mathcal{M}_{+}({{\mathbb{R}}_{\ast}^{d}})$ endowed with the weak-$\ast$
topology is called a weak solution to \eqref{B4} with initial value
$f_{0}$ if $f(0,\cdot)=f_{0}(\cdot)$ and for each $1<T<\infty$%
\begin{equation}
\sup_{t\in\left[  1/T,T\right]  }[\int_{\left\{  \left\vert x\right\vert
\geq1\right\}  }|x|^{\gamma+p}f(dx,t)+\int_{\left\{  \left\vert x\right\vert
\leq1\right\}  }|x|^{1-p}f(dx,t)]<\infty, \label{eq:bound_gamma_moment}%
\end{equation}
\begin{equation}
\int_{{\mathbb{R}}_{\ast}^{d}}xf(dx,t)=\int_{{\mathbb{R}}_{\ast}^{d}}%
xf_{0}(dx),\quad t>0 \,, \label{eq:conservation}%
\end{equation}
and, for all test functions $\varphi\in C_{c}^{1}({{\mathbb{R}}_{\ast}%
^{d}\times(0,\infty)})$
the following identity holds
\begin{align}
0  &  =\int_{0}^{\infty}\int_{{{\mathbb{R}}_{\ast}^{d}}}f(dx,t)\partial
_{t}\varphi(x,t) dt\nonumber\\
&  +\frac{1}{2}\int_{0}^{\infty}\int_{{{\mathbb{R}}_{\ast}^{d}}}%
\int_{{{\mathbb{R}}_{\ast}^{d}}}K(x,y)f(dx,t)f(dy,t)[\varphi(x+y,t)-\varphi
(x,t)-\varphi(y,t)]dt. \label{eq:weak}%
\end{align}

\end{definition}

\begin{remark}
Notice that the condition \eqref{eq:bound_gamma_moment} is equivalent to
\begin{equation}
\label{eq:a_bound_gamma_moment}\sup_{t\in\left[  1/T,T\right]  }%
[\int_{\left\{  \left\vert x\right\vert \geq a\right\}  }|x|^{\gamma
+p}f(dx,t)+\int_{\left\{  \left\vert x\right\vert \leq a\right\}  }%
|x|^{1-p}f(dx,t)]<\infty,
\end{equation}
for any $a>0$. 
 In order to define the solutions it would be enough to impose an integrability condition in $t$ and $x$. However we decided to stick to the stronger condition \eqref{eq:bound_gamma_moment} as it allows us to use estimates derived in \cite{EM06}. 
\end{remark}

Note that in Definition \ref{def:time-dep_sol} we allow only solutions with
finite mass that is conserved over time for each component due to condition
\eqref{eq:conservation}. The assumption \eqref{eq:bound_gamma_moment} ensures
that all the integrals appearing in \eqref{eq:weak} are well-defined for
kernels satisfying the upper bound \eqref{eq:condK_sym2}. Indeed, the last
term in (\ref{eq:weak}) can be estimated by splitting the domain of
integration into two regions defined by $\{|y| \leq|x| \}$ and $\{|x| > |y|
\}$. Using a symmetrization argument, the integral over the second region can
be estimated by the integral over the first region and therefore,
\eqref{eq:weak} can be estimated as
\begin{equation}
\int_{0}^{\infty}\int_{\left\{  \left\vert y\right\vert \leq\left\vert
x\right\vert \right\}  }\int K(x,y)f(dx,t)f(dy,t)\left\vert \varphi
(x+y,t)-\varphi(x,t)-\varphi(y,t)\right\vert dt. \label{A1}%
\end{equation}

Since $\varphi\in C_{c}^{1}({{\mathbb{R}}_{\ast}^{d}\times}\left(
0,\infty\right)  )$, there exists a function $\psi\in C_{c}({{\mathbb{R}%
}_{\ast}^{d}})$ such that $|\varphi(x+y,t)-\varphi(x,t)|\leq\psi(x)|y|$. Let
$\supp\psi,\supp\varphi\subset\{x\ |\ \frac{1}{L}<|x|<L\}$, for some positive
constant $L>1$. Using the upper bound \eqref{eq:condK_sym2} for the kernel
$K$, the term in \eqref{eq:weak} involving $\varphi(x+y,t)-\varphi
(x,t)-\varphi(y,t)$ can be estimated by%
\begin{align*}
&  \iint_{\left\{  \left\vert y\right\vert \leq\left\vert x\right\vert
\right\}  }K(x,y)f(dx,t)f(dy,t)|\varphi(x+y,t)-\varphi(x,t)| \\
& \quad \quad +\iint_{\left\{
\left\vert y\right\vert \leq\left\vert x\right\vert \right\}  }%
K(x,y)f(dx,t)f(dy,t)\left\vert \varphi(y,t)\right\vert = \\
&  =\int_{{\{x\in{\mathbb{R}}_{\ast}^{d}\ :\ }\frac{1}{L}{\leq|x|\leq L\}}%
}\int_{\left\{  \left\vert y\right\vert \leq\left\vert x\right\vert \right\}
}\left(  \cdot\cdot\cdot\right)  +\int_{{\{x\in{\mathbb{R}}_{\ast}^{d}%
\ |\ }\frac{1}{L}{\leq|y|\leq L\}}}\int_{\left\{  \left\vert y\right\vert
\leq\left\vert x\right\vert \right\}  }\left(  \cdot\cdot\cdot\right) \\
&  \leq\int_{{\{x\in{\mathbb{R}}_{\ast}^{d}\ :\ }\frac{1}{L}{\leq|x|\leq L\}}%
}\int_{\left\{  \left\vert y\right\vert \leq\left\vert x\right\vert \right\}
}K(x,y)\psi(x)|y|f(dx,t)f(dy,t)+\\
&  +\int_{{\{x\in{\mathbb{R}}_{\ast}^{d}\ |\ }\frac{1}{L}{\leq|y|\leq L\}}%
}\int_{\left\{  \frac{1}{L}\leq\left\vert x\right\vert \right\}
}K(x,y)f(dx,t)f(dy,t)\\
&  \leq C\int_{\{y\in{{\mathbb{R}}_{\ast}^{d}}\ |\ |y|\leq L\}}\left(
|y|^{-p}|y|+|y|^{\gamma+p}|y|\right)  f(dy,t)+C\int_{\left\{  \frac{1}{L}%
\leq\left\vert x\right\vert \right\}  }\left(  |x|^{\gamma+p}+|x|^{-p}\right)
f(dx,t)\\
&  \leq C\int_{\{y\in{{\mathbb{R}}_{\ast}^{d}}\ |\ |y|\leq L\}}|y|^{1-p}%
f(dy,t)+C\int_{\left\{  \frac{1}{L}\leq\left\vert x\right\vert \right\}
}|x|^{\gamma+p}f(dx,t)<\infty \,.
\end{align*}
The finiteness of the integrals follows from the assumption
\eqref{eq:bound_gamma_moment} (more precisely it follows from
\eqref{eq:a_bound_gamma_moment} with $a=1/L, a=L$). Notice that the constants
$C$ depend on $L.$

We remark that using a standard limit argument, we obtain that (\ref{eq:weak})
implies that the following identity holds for test functions $\varphi\in
C_{c}^{1}({{\mathbb{R}}_{\ast}^{d}\times\lbrack0,\infty)})$,%
\begin{align}
0  &  =\int_{{{\mathbb{R}}_{\ast}^{d}}}f_{0}(dx)\varphi(x,0)+\int_{0}^{\infty
}\int_{{{\mathbb{R}}_{\ast}^{d}}}f(dx,t)\partial_{t}\varphi(x,t)dt \nonumber\\
&  +\frac{1}{2}\int_{0}^{\infty}\int_{{{\mathbb{R}}_{\ast}^{d}}}%
\int_{{{\mathbb{R}}_{\ast}^{d}}}K(x,y)f(dx,t)f(dy,t)[\varphi(x+y,t)-\varphi
(x,t)-\varphi(y,t)]dt. \label{eq:weakbis}%
\end{align}

The localization result in Theorem \ref{TheorSelfSimLoc} concerns self-similar
solutions. In the next definition we collect the properties required for the
self similar solutions, and more specifically for the self-similar profiles
$F$ (cf.\ (\ref{MultSelfSimForm})) which we need in order to derive the
localization result.

\begin{definition}
\label{def:self_sim_sol} Let the kernel $K$ satisfy \eqref{eq:condK_cont}, (\ref{eq:condK_sym1}),
(\ref{eq:condK_sym2}), (\ref{eq:Phi_kernel}), (\ref{KernIneq}) as well as the homogeneity condition (\ref{eq:condK_homogeneity}). We say that a measure
$F\in\mathcal{M}_{+}({{\mathbb{R}}_{\ast}^{d}})$ satisfying\
\begin{equation}
\int_{\left\{  \left\vert \xi\right\vert \geq1\right\}  }|\xi|^{\gamma
+p}F(d\xi)+\int_{\left\{  \left\vert \xi\right\vert <1\right\}  }|\xi
|^{1-p}F(d\xi)+\int_{{{\mathbb{R}}_{\ast}^{d}}}|\xi|F(d\xi)<\infty
\label{eq:bound_gamma_moment_selfsim}%
\end{equation}
is a self-similar profile to \eqref{B4} if the following identity
holds
\begin{align}
& 0=  \frac{1}{2}\int_{{{\mathbb{R}}_{\ast}^{d}}}\int_{{{\mathbb{R}}_{\ast}^{d}}%
}F(d\xi)F(d\eta)K(\xi,\eta)[\psi(\xi+\eta)-\psi(\xi)-\psi(\eta
)]\nonumber\\
&  +\frac{1}{1-\gamma}\int_{{{\mathbb{R}}_{\ast}^{d}}}F(d\xi)\left[  \psi
-\xi\cdot\partial_{\xi}\psi\right]  \label{WeakSelfSim}
\end{align}
for all $\psi\in C_{c}^{1}({{\mathbb{R}}_{\ast}^{d}})$.
\end{definition}

The finiteness of the integrals in
(\ref{eq:bound_gamma_moment_selfsim}) provide the integrability required to
ensure that the integrals in (\ref{WeakSelfSim}) are well defined. On the
other hand, the finiteness of the last integral in
(\ref{eq:bound_gamma_moment_selfsim}) implies that the total number of
monomers associated to the function $f$ defined by means of
(\ref{MultSelfSimForm}) is finite. Notice that a self-similar profile can be
interpreted as a weak solution of the equation%
\begin{equation}
\mathbb{K}[F]+\frac{1}{1-\gamma}\xi\cdot\partial_{\xi}F+\frac{d+1}{1-\gamma
}F=0\ \ ,\ \ \xi\in{{\mathbb{R}}_{\ast}^{d}} \label{SelfSimK}%
\end{equation}
with $\mathbb{K}[F]$ as in (\ref{coagKern}). Equation (\ref{SelfSimK}) may be
obtained formally from the coagulation equation \eqref{B4} using the
change of variables \eqref{MultSelfSimForm}. This change of variables can be made precise by noticing that  given a self-similar profile $F$ in the sense of Definition
\ref{def:self_sim_sol} we can obtain a weak solution $f$ of (\ref{B4}) in
the sense of Definition \ref{def:time-dep_sol} by requiring
\begin{equation}\label{MultSelfSimFormRigorous}
\iint_{\R_*^d \times [0,\infty)} f(x,t)\varphi(x,t)dxdt = \iint_{\R_*^d \times [0,\infty)} \varepsilon_t F\left(  \xi\right) \varphi(\xi (\varepsilon_t)^{-1} ,t)d\xi dt,
\end{equation}
  for any $\varphi \in C_c(\R_*^d \times [0,\infty))$ and where $\varepsilon_t = (1+t)^{-\frac{1}{1-\gamma}}$. We have the following result.

\begin{proposition}
\label{SelfSimProfSol} Suppose that $K$ satisfies
(\ref{eq:condK_cont}), (\ref{eq:condK_sym1}), (\ref{eq:condK_sym2}), (\ref{eq:Phi_kernel}),
(\ref{KernIneq}) as well as the homogeneity condition (\ref{eq:condK_homogeneity}). Let
us assume also that $F\in\mathcal{M}_{+}({{\mathbb{R}}_{\ast}^{d}})$ is a
self-similar profile in the sense of Definition \ref{def:self_sim_sol}. We
define $f$ for $t\geq0$ as in (\ref{MultSelfSimForm}) (cf.\ also \eqref{MultSelfSimFormRigorous}).  Then, $f$ is a weak
solution of (\ref{B4}) in the sense of Definition
\ref{def:time-dep_sol} with initial value $f_{0}=F$, satisfying the moment bounds  for all $T>0$  
\begin{equation*}
\sup_{t\in\left[ 0,T\right]  }[\int_{\left\{  \left\vert x\right\vert
\geq1\right\}  }|x|^{\gamma+p}f(dx,t)+\int_{\left\{  \left\vert x\right\vert
\leq1\right\}  }|x|^{1-p}f(dx,t)]<\infty.
\end{equation*}
Moreover, $f$ is
invariant under the following group of transformations:
\begin{equation}\label{eq:f_lambda}
f_{\lambda}\left(  x,t\right)  =\lambda^{d+1}f\left(  \lambda x,\lambda
^{1-\gamma}\left(  t+1\right)  -1\right)  \ \ ,\ \ \ \lambda>0.
\end{equation}
\end{proposition}

\begin{remark}
We notice that $f_\lambda$ and $f$ are measures and 
 \eqref{eq:f_lambda} should be interpreted as follows, for any $\varphi \in C_c(\R_*^d \times [0,\infty))$,
$$\iint_{\R_*^d \times [0,\infty)} f_\lambda(x,t)\varphi(x,t)dxdt = \lambda^{\gamma} \iint_{\R_*^d \times [0,\infty)} f(x,t)\varphi \left(\frac{x}{\lambda}, \frac{t+1}{\lambda^{1-\gamma}}-1\right) dxdt.$$
\end{remark}

\begin{proof}
Given that $f$ is defined in (\ref{MultSelfSimForm}) it follows from
(\ref{eq:bound_gamma_moment_selfsim}) that (\ref{eq:bound_gamma_moment})
holds. We now compute the right-hand side of (\ref{eq:weakbis}). Notice that
(\ref{MultSelfSimForm}) implies that $f_{0}\left(  x\right)  =F\left(
x\right)  .$ Then%
\begin{align*}
&  \int_{{{\mathbb{R}}_{\ast}^{d}}}f_{0}(dx)\varphi(x,0)+\int_{0}^{\infty}%
\int_{{{\mathbb{R}}_{\ast}^{d}}}f(dx,t)\partial_{t}\varphi(x,t)\\
&  +\frac{1}{2}\int_{0}^{\infty}\int_{{{\mathbb{R}}_{\ast}^{d}}}%
\int_{{{\mathbb{R}}_{\ast}^{d}}}K(x,y)f(dx,t)f(dy,t)[\varphi(x+y,t)-\varphi
(x,t)-\varphi(y,t)]dt\\
&  =\int_{{{\mathbb{R}}_{\ast}^{d}}}F(dx)\varphi(x,0)+\int_{0}^{\infty}%
\int_{{{\mathbb{R}}_{\ast}^{d}}}f(dx,t)\partial_{t}\varphi(x,t)\\
&  +\frac{1}{2}\int_{0}^{\infty}\int_{{{\mathbb{R}}_{\ast}^{d}}}%
\int_{{{\mathbb{R}}_{\ast}^{d}}}K(x,y)f(dx,t)f(dy,t)[\varphi(x+y,t)-\varphi
(x,t)-\varphi(y,t)]dt\\
&  :=J .
\end{align*}

Given a test function $\varphi\in C_{c}^{1}({{\mathbb{R}}_{\ast}^{d}%
\times\left[  0,\infty\right)  })$ we define $\psi\left(  \xi,\tau\right)  $
by means of
\begin{equation}
\varphi\left(  x,t\right)  =(\varepsilon_t)^{-1}
\psi\left(  \xi,\tau\right)  \label{eq:change_phi}
\end{equation}
where%
\[
\xi=x \varepsilon_t\ \ ,\ \ \tau
=\log\left(  t+1\right)\ \ ,\ \ \varepsilon_t = (1+t)^{-\frac{1}{1-\gamma}}   .
\]
Then, using also (\ref{MultSelfSimForm}) as well as   the homogeneity condition
(\ref{eq:condK_homogeneity}) and $ d\tau  = \frac{dt}{t+1}  $ we obtain 
\begin{align*}
J  &  =\int_{{{\mathbb{R}}_{\ast}^{d}}}F(d\xi)\psi(\xi,0)+\int_{0}^{\infty
}d\tau\int_{{{\mathbb{R}}_{\ast}^{d}}}F(d\xi)\left[  \partial_{\tau}\psi
+\frac{\psi}{1-\gamma}-\frac{1}{1-\gamma}\xi\cdot\partial_{\xi}\psi\right]
(\xi,\tau)\\
&  +\frac{1}{2}\int_{0}^{\infty}d\tau\int_{{{\mathbb{R}}_{\ast}^{d}}}%
\int_{{{\mathbb{R}}_{\ast}^{d}}}K\left(  \xi,\eta\right)  F(d\xi)F(d\eta
)[\psi(\xi+\eta,\tau)-\psi(\xi,\tau)-\psi(\eta,\tau)].
\end{align*}
Employing (\ref{WeakSelfSim}), we then find
\[
J=\int_{{{\mathbb{R}}_{\ast}^{d}}}F(d\xi)\psi(\xi,0)+\int_{0}^{\infty}%
d\tau\int_{{{\mathbb{R}}_{\ast}^{d}}}F(d\xi)\partial_{\tau}\psi=\int
_{{{\mathbb{R}}_{\ast}^{d}}}F(d\xi)\psi(\xi,0)-\int_{{{\mathbb{R}}_{\ast}^{d}%
}}F(d\xi)\psi(\xi,0)=0
\]
and the result follows.
\end{proof}

\section{Proof of the localization results (Theorems \ref{mainRoughImproved},
\ref{mainRough}, \ref{TheorSelfSimLoc})}

\subsection{Mass localization along a ray in time-dependent solutions}

In this subsection we prove Theorems \ref{mainRoughImproved} and
\ref{mainRough}. To this end, it is convenient to rewrite the function $f$ in
the theorems using the set of self-similar variables 
\begin{equation}
f\left(  x,t\right)  =(\varepsilon_t)^{1+d}F\left(  \xi,\tau\right)  ,\ \ \xi=x\varepsilon_t,\ \ \ \tau=\log(t+1), \ \ \varepsilon_t = (t+1)^{-\frac{1}{1-\gamma}} \label{MultSelfSimFormT}%
\end{equation}
where $f$ is a weak solution of \eqref{B4} in the sense of definition
\ref{def:time-dep_sol}. In order to prove the localization results it will be
convenient to further define a set of simplicial coordinates
\begin{equation}
\left(  \rho,\theta\right)  =\left(  \left\vert \xi\right\vert ,\frac{\xi
}{\left\vert \xi\right\vert }\right)  \ \ \text{with }\rho\in\mathbb{R}%
_{*}\ \ ,\ \ \theta\in\Delta^{d-1} \label{rhoTheta}%
\end{equation}
where $\Delta^{d-1}$ is as in (\ref{eq:simplex}). A similar system of
coordinates has been used also in \cite{FLNV2}. We denote as $d\nu\left(
\theta\right)  $ the $\left(  d-1\right)  $ dimensional Hausdorff measure
restricted to $\Delta^{d-1}$. 
We thus have that $d\xi
=\frac{\rho^{d-1}}{\sqrt{d}}d\rho d\nu\left(  \theta\right)$ (cf.~\cite{FLNV2}).  We can then
define measures $G(\tau)\in\mathcal{M}_{+}\left(  \mathbb{R}_{*}\times\Delta
^{d-1}\right)  $ by requiring that 
\begin{equation}
\int_{ \mathbb{R}_{*}\times\Delta^{d-1}} \psi(\rho,\theta) 
G(\rho,\theta,\tau)\rho^{d-1} d\rho d\tau\left(  \theta\right) 
=\int_{\R^d_\ast} \psi(|\xi|,\xi/|\xi|) F(\xi,\tau) d\xi\,, 
\label{PolDist}
\end{equation}
for all test functions $\psi\in  C_c(\R_*\times \Delta^{d-1})$ and
with $F$ defined by \eqref{MultSelfSimFormT}. Notice that in the case in which
$F\left(  \cdot,\tau\right)  $ is absolutely continuous this implies that
$G\left(  \cdot, \cdot,\tau\right)  $ is also absolutely continuous and we
have the following relation between the corresponding densities%
\[
F\left( \xi,\tau\right)  =\sqrt{d}\,G\left(|\xi|,\xi/|\xi|,\tau\right)\,.
\]

We then use (\ref{PolDist}) and the property $t+1 = e^{\tau}$ and follow an argument similar to that in the proof of Proposition
\ref{SelfSimProfSol}.  This allows us to rewrite
(\ref{eq:weak}) as%
\begin{align}
&  0=\int_{0}^{\infty}d\tau\int_{ \mathbb{R}_{*}\times\Delta^{d-1}}G\left(
\rho,\theta,\tau\right)  \left[  \partial_{\tau}\tilde{\psi}+\frac{\tilde
{\psi}}{1-\gamma}-\frac{1}{1-\gamma}\rho\partial_{\rho}\tilde{\psi}\right]
(\rho,\theta,\tau)d\Omega\nonumber\\
&  +\frac{1}{2}\int_{0}^{\infty}d\tau\int_{ \mathbb{R}_{*}\times\Delta^{d-1}}\int_{ \mathbb{R}_{*}\times\Delta^{d-1}}\tilde{K}(\rho,\theta,r,\sigma
,\tau)G\left(  \rho,\theta,\tau\right)  G\left(  r,\sigma,\tau\right)
\Psi(\rho,r,\theta,\sigma,\tau)d\Omega d\tilde{\Omega} \label{WeakFormSimp}%
\end{align}
with
\begin{equation}
\tilde{K}(\rho,\theta,r,\sigma,\tau)=e^{-\frac{\gamma}{1-\gamma}\tau
}K(e^{\frac{1}{1-\gamma}\tau}\rho\theta,e^{\frac{1}{1-\gamma}\tau}r\sigma)
\label{Ktilde}%
\end{equation}
\[
d\Omega=\rho^{d-1}d\rho d\nu(\theta),\quad d\tilde{\Omega}=r^{d-1}%
drd\nu(\sigma),
\]%
\begin{equation}
\Psi(\rho,r,\theta,\sigma,\tau)=\tilde{\psi}\left(  \rho+r,\frac{\rho}{\rho
+r}\theta+\frac{r}{\rho+r}\sigma,\tau\right)  -\tilde{\psi}(\rho,\theta
,\tau)-\tilde{\psi}(r,\sigma,\tau) \label{PsiTest}%
\end{equation}
where we write $\tilde{\psi}(\rho,\theta,\tau)=\psi(\xi,t)$ with $\psi$
defined in \eqref{eq:change_phi}$,\ \xi=\rho\theta,\ \eta=r\sigma$ and \\
$G_{0}\left(  \rho,\theta\right)  \rho^{d-1}d\rho d\nu\left(  \theta\right)
=f_{0}(dx)$ for each $\theta,\ \sigma\in\Delta^{d-1}.$

Notice that (\ref{eq:condK_sym1}), (\ref{eq:condK_sym2}) imply the estimate%
\begin{equation}
c_{1}(\rho+r)^{\gamma}\Phi_{p}\left(  \frac{\rho}{\rho+r}\right)  \leq\tilde
K(\rho,\theta,r,\sigma,\tau)\leq c_{2}(\rho+r)^{\gamma}\Phi_{p}\left(
\frac{\rho}{\rho+r}\right)  \label{KestimSimplex}%
\end{equation}
where $\Phi_{p}$ is as in (\ref{eq:Phi_kernel}).

The following Lemma, which will be used to obtain Theorem \ref{mainRough}, has
been proved in \cite{FLNV2}. For this reason, we will just state the result
and refer to \cite{FLNV2} for the proof.

\begin{lemma}
\label{LemmSplit} There is a constant $C_{d}>0$ which depends only on the
dimension $d\geq1$ such that, for any probability measure $\lambda
\in\mathcal{M}_{+}\left(  \Delta^{d-1}\right)  $ and any pair of parameters
$\varepsilon,\delta\in(0,1)$ at least one of the following alternatives holds true:

\begin{itemize}
\item[(i)] There exists a measurable set $A\subset\Delta^{d-1}$ with
$\mathrm{diam}\left(  A\right)  \leq\varepsilon$ such that $\int_{A}%
\lambda(d\theta)>1-\delta$.

\item[(ii)] $\int_{\Delta^{d-1}}\lambda\left(  d\theta\right)  \int
_{\Delta^{d-1}}\lambda\left(  d\sigma\right)  \left\Vert \theta-\sigma
\right\Vert ^{2}\geq C_{d}\delta\varepsilon^{d+1}$ where $\left\Vert
\cdot\right\Vert $ is the Euclidean distance.
\end{itemize}
\end{lemma}

A corollary of Lemma \ref{LemmSplit} that will be used in the proof of Theorem
\ref{TheorSelfSimLoc} is the following result.

\begin{lemma}
\label{LemmSplitDelta} Suppose that $\lambda\in\mathcal{M}_{+}\left(
\Delta^{d-1}\right)  $ is a probability measure such that%
\begin{equation}
\int_{\Delta^{d-1}}\lambda\left(  d\theta\right)  \int_{\Delta^{d-1}}%
\lambda\left(  d\sigma\right)  \left\Vert \theta-\sigma\right\Vert ^{2}=0.
\label{LambIdent}%
\end{equation}

Then, there exists $\bar{\theta}\in\Delta^{d-1}$ such that%
\begin{equation}
\lambda= \delta_{\bar{\theta}}\,.
\label{LambDelt}%
\end{equation}

\end{lemma}

\begin{proof}
We apply Lemma \ref{LemmSplit} for a sequence of values $\varepsilon
_{n}=\delta_{n}=2^{-n}$ with $n\in\mathbb{N}$. Due to (\ref{LambIdent}) we
have that the alternative (ii) in Lemma \ref{LemmSplit} does not take place.
Therefore (i) holds, and
for each $n$, we can pick a point $\theta_n$ from the corresponding set $A$.  By compactness of $\Delta^{d-1}$, we can find a convergent subsequence
with a limit point $\bar{\theta}\in\Delta^{d-1}$,
and it can be checked that then also (\ref{LambDelt}) holds.
\end{proof}

We can now prove Theorem \ref{mainRough}.

\begin{proofof}
[Proof of Theorem \ref{mainRough}]The conservation of the total number of
monomers (\ref{eq:conservation}), combined with the definition of $F$ in
(\ref{MultSelfSimFormT}) and the definition of $G$ in (\ref{PolDist}) implies
that
\begin{equation}
\int_{{\mathbb{R}}_{*}\times\Delta^{d-1}}\rho G(\rho,\theta,\tau)d\Omega
=m_{0} := \int_{{\mathbb{R}}_{*}\times\Delta^{d-1}}\rho G_{0}\left(  \rho
,\theta\right)  d\Omega>0 \,.\label{eq:mass1}%
\end{equation}
We recall that the moment estimate \eqref{eq:moment_estimate} holds for $t\geq 1$. Since $ \tau=\log(t+1)$  we will assume $\tau \geq \ln 2$ throughout the proof.

On the other hand, the assumption (\ref{eq:moment_estimate}) with
(\ref{MultSelfSimFormT}), (\ref{PolDist}) yields
\begin{equation}
\int_{{\mathbb{R}}_{*}\times\Delta^{d-1}}\rho^{k}G(\rho,\theta,\tau
)d\Omega\leq C,\quad k\in\lbrack1/a,a]\ \ \text{for some }%
a>1,
\label{eq:sec_moment}%
\end{equation}
where $C>0.$
In addition, using again (\ref{MultSelfSimFormT}), (\ref{PolDist}) as well as
the estimate (\ref{eq:bound_gamma_moment}) we obtain%
\begin{equation}
\int_{\left[  {\mathbb{R}}_{*}\times\Delta^{d-1}\right]  \cap\left\{  \rho
\geq1\right\}  }\rho^{\gamma+p}G(\rho,\theta,\tau)d\Omega+\int_{\left[
{\mathbb{R}}_{*}\times\Delta^{d-1}\right]  \cap\left\{  \rho\leq1\right\}
}\rho^{1-p}G(\rho,\theta,\tau)d\Omega\leq C\left(  T\right)  ,\ \ \ln 2\leq\tau\leq T \label{MomEstG}%
\end{equation}
for any given $T>1,$ with $C(T)$ a constant that depends on $T$.

In the definition of weak solution (cf.\ Definition \ref{def:time-dep_sol}) we
have assumed that the mass vector is conserved. Using the measure $G$ we then
obtain the following form of the conservation of mass%
\begin{equation}
\frac{1}{m_{0}}\int_{{\mathbb{R}}_{*}\times\Delta^{d-1}}G(\rho,\theta
,\tau)\rho\theta d\Omega=\theta_{0}, \label{eq:mass_conservation}%
\end{equation}
with%
\begin{equation}
\theta_{0}:=\frac{1}{m_{0}}\int_{{\mathbb{R}}_{*}\times\Delta^{d-1}}G_{0}%
(\rho,\theta)\rho\theta d\Omega. \label{eq:def_theta0}%
\end{equation}

Using (\ref{eq:sec_moment}) and (\ref{MomEstG}) we can readily see that all
the terms appearing in (\ref{WeakFormSimp}) are well defined for any
$\tilde{\psi}\in C_{c}^{1}\left(  \mathbb{R}_{*}\times\Delta^{d-1}%
\times\left( 0,\infty\right)  \right)$. Moreover, using an approximation
argument as well as (\ref{eq:sec_moment}) and (\ref{MomEstG}) it follows that
(\ref{WeakFormSimp}) holds true for any test function $\tilde{\psi}\in
C^{1}\left(  \mathbb{R}_{*}\times\Delta^{d-1}\times\left[  \ln(2),\infty\right)
\right)  $ whose support is contained in $\left\{  \tau:\tau\in\left[
\ln(2),\tau^{\ast}\right]  \right\}  $, for some $\tau^{\ast}\in\left(
\ln(2),\infty\right) $, and satisfying
\begin{equation}
\left\vert \tilde{\psi}(\rho,\theta,\tau)\right\vert +\rho\left\vert
\partial_\rho \tilde{\psi}(\rho,\theta,\tau)\right\vert
+\left\vert \partial_{\tau}\tilde{\psi}(\rho,\theta,\tau)\right\vert
+\left\vert \nabla_{\theta}\tilde{\psi}(\rho,\theta,\tau)\right\vert \leq
C\rho\ \ , \label{TestLinEstim}%
\end{equation}
for any $ \rho\in\mathbb{R}_{*},\ \theta\in\Delta^{d-1}\text{ and }\tau
\in\left[  \ln(2),\tau^{\ast}\right] $, and for some $C>0.$ Indeed, we can consider a sequence of compactly supported test
functions $\tilde{\psi}_{n}(\rho,\theta,\tau)=\zeta_{n}\left(  \rho\right)
\tilde{\psi}(\rho,\theta,\tau)$ with $\tilde{\psi}$\ satisfying
(\ref{TestLinEstim}), $\zeta_{n}\in C^{\infty}\left(  0,\infty\right)  $  such
that
\[
\zeta_{n}\left(  \rho\right)  =\left\{
\begin{array}
[l]{cc}%
0 & \text{ \quad for \ }0<\rho\leq\frac{1}{2n}\text{\quad and \quad}\rho
\geq2n\\
1 & \text{ \quad for \quad}\rho\in\left[  \frac{1}{n},n\right]
\end{array}
\right.
\]
and $0\leq\zeta_{n}{}^{\prime}\left(  \rho\right)  \leq4n$ for $\rho\in\left[
\frac{1}{2n},\frac{1}{n}\right]  $ and $0\leq-\zeta_{n}{}^{\prime}\left(
\rho\right)  \leq\frac{2}{n}$ for $\rho\in\left[  n,2n\right]  .$ By
assumption the identity (\ref{WeakFormSimp}) holds with the test functions
$\tilde{\psi}_{n}.$ We consider the limit as $n\rightarrow\infty$ of this
sequence of identities. Notice that (\ref{TestLinEstim}) combined with the
properties of $\zeta_{n}$ imply the estimate%
\[
\left\vert \partial_{\tau}\tilde{\psi}_{n}+\frac{\tilde{\psi}_{n}}{1-\gamma
}-\frac{1}{1-\gamma}\rho\partial_{\rho}\tilde{\psi}_{n}\right\vert \leq C\rho
\]
with $C$ independent of $n.$ We can then take the limit as $n\rightarrow
\infty$ of the first term on the right-hand side of (\ref{WeakFormSimp}) (with
$\tilde{\psi}$ replaced by $\tilde{\psi}_{n}$) due to Lebesgue's dominated
convergence Theorem as well as (\ref{eq:mass_conservation}). On the other
hand, using (\ref{TestLinEstim})
and the properties of $\zeta_{n},$ we can estimate the functions $\Psi_{n}(\rho,r,\theta,\sigma,\tau)$ that
are defined using $\tilde{\psi}=\tilde{\psi}_{n}$ as
\begin{equation}
\label{eq:Psin}\left\vert \Psi_{n}(\rho,r,\theta,\sigma,\tau)\right\vert \leq
C\min\left\{  \rho,r\right\}
\end{equation}
where $C$ is independent of $n.$ To see this we can restrict ourselves to the
case in which $r\leq\rho.$\ We can then estimate the term $\tilde{\psi}%
_{n}\left(  r,\sigma,\tau\right)  $ in the formula of $\Psi_{n}$ (cf.
(\ref{PsiTest})) by $Cr$\ and the difference $\tilde{\psi}\left(  \rho
+r,\frac{\rho}{\rho+r}\theta+\frac{r}{\rho+r}\sigma,\tau\right)  -\tilde{\psi
}(\rho,\theta,\tau)$ can be estimated also as $Cr$ using Taylor's formula and
(\ref{TestLinEstim}). Using (\ref{KestimSimplex}) and (\ref{eq:Phi_kernel}) it
then follows that the integrands in the second term of (\ref{WeakFormSimp})
(with $\tilde{\psi}=\tilde{\psi}_{n}$) can be estimated by an integrable
function independent of $n,$ due to (\ref{MomEstG}) and
(\ref{eq:mass_conservation}). Indeed, the previous estimate \eqref{eq:Psin}
yields the bound%
\begin{align}
&  \tilde{K}(\rho,\theta,r,\sigma,\tau)G\left(  \rho,\theta,\tau\right)
G\left(  r,\sigma,\tau\right)  \Psi_{n}(\rho,r,\theta,\sigma,\tau)\nonumber\\
&  \leq C\left(  \rho^{\gamma+p}r^{-p}+r^{\gamma+p}\rho^{-p}\right)  rG\left(
\rho,\theta,\tau\right)  G\left(  r,\sigma,\tau\right) \nonumber\\
&  \leq C\left(  \rho^{\gamma+p}r^{1-p}+\rho^{-p}r^{\gamma+p+1}\right)
G\left(  \rho,\theta,\tau\right)  G\left(  r,\sigma,\tau\right)
\label{KernLebDom}%
\end{align}
where in the last inequality we use that $r\leq\rho.$ This inequality gives an
estimate for the integrand of the last term in (\ref{WeakFormSimp}) by means
of an integrable function for $0<r\leq\rho\leq1,$ due to the assumption
(\ref{eq:bound_gamma_moment}). In order to estimate the regions where
$0<r\leq1\leq\rho$ we use the fact that $\gamma+p<1$ and $p>-1$ to obtain the
following estimate for this range of values of $r$ and $\rho$%
\begin{align*}
&  \tilde{K}(\rho,\theta,r,\sigma,\tau)G\left(  \rho,\theta,\tau\right)
G\left(  r,\sigma,\tau\right)  \Psi_{n}(\rho,r,\theta,\sigma,\tau)\\
&  \leq C\left(  \rho r^{1-p}+\rho r^{\gamma+p+1}\right)  G\left(  \rho
,\theta,\tau\right)  G\left(  r,\sigma,\tau\right) .
\end{align*}

We can then combine the conservation of mass estimate and
(\ref{eq:bound_gamma_moment}) to obtain that the integrand is bounded by an
uniformly integrable function independent of $n$ for this range of parameters.
It remains to estimate the range of values $1\leq r\leq\rho.$ We distinguish
the two cases $p\geq0$ and $-1<p<0.$ In the first case, we use
(\ref{KernLebDom}) to obtain, using also that $\gamma+p < 1$ and $\gamma<1$
(cf.\ (\ref{KernIneq}))%
\begin{align*}
&  \tilde{K}(\rho,\theta,r,\sigma,\tau)G\left(  \rho,\theta,\tau\right)
G\left(  r,\sigma,\tau\right)  \Psi_{n}(\rho,r,\theta,\sigma,\tau)\\
&  \leq C\left(  \rho r+\rho^{-p}r^{\gamma+p+1}\right)  G\left(  \rho
,\theta,\tau\right)  G\left(  r,\sigma,\tau\right) \\
&  \leq C\left(  \rho r+\rho r^{\gamma}\right)  G\left(  \rho,\theta
,\tau\right)  G\left(  r,\sigma,\tau\right)
\end{align*}
where in the last inequality we used that $r^{p+1}\leq\rho^{p+1}.$ In the case
$-1< p<0$ we obtain%
\begin{align*}
&  \tilde{K}(\rho,\theta,r,\sigma,\tau)G\left(  \rho,\theta,\tau\right)
G\left(  r,\sigma,\tau\right)  \Psi_{n}(\rho,r,\theta,\sigma,\tau)\\
&  \leq C\left(  \rho^{\gamma}\rho^{-\left\vert p\right\vert }r^{\left\vert
p\right\vert }r+r^{\gamma-\left\vert p\right\vert }\rho^{\left\vert
p\right\vert }r\right)  G\left(  \rho,\theta,\tau\right)  G\left(
r,\sigma,\tau\right) \\
&  \leq C\left(  \rho^{\gamma}r+r^{\gamma-\left\vert p\right\vert }%
\rho^{\left\vert p\right\vert }\rho^{1-\left\vert p\right\vert }r^{\left\vert
p\right\vert }\right)  G\left(  \rho,\theta,\tau\right)  G\left(
r,\sigma,\tau\right) \\
&  \leq C\left(  \rho^{\gamma}r+r^{\gamma}\rho\right)  G\left(  \rho
,\theta,\tau\right)  G\left(  r,\sigma,\tau\right) .
\end{align*}

The right-hand side of this inequality is uniformly integrable for $1\leq
r\leq\rho$ since $\gamma<1.$

We can then take the limit in both terms of (\ref{WeakFormSimp}) (with
$\tilde{\psi}=\tilde{\psi}_{n}$). It then follows that (\ref{WeakFormSimp})
holds for any function $\tilde{\psi}$ satisfying (\ref{TestLinEstim}).

We now approximate a test function with the form $\tilde{\psi}(\rho
,\theta,\tau)\chi_{\left[  \bar{\tau}_{1},\bar{\tau}_{2}\right]  }$ with
$\bar{\tau}_{2}\geq\bar{\tau}_{1} \geq \ln(2) $ and $\tilde{\psi}\in C^{1}\left(\mathbb{R}_{\ast}\times\Delta^{d-1}\times\left[  \ln(2),\infty\right)  \right)$ 
by
means of a sequence of functions $\tilde{\psi}_{n}(\rho,\theta,\tau
)=\tilde{\psi}(\rho,\theta,\tau)\zeta_{n}\left(  \tau\right)  $ with
$\zeta_{n}\in C_{c}^{1}\left(  \left(  0,\infty\right)  \right)  ,$ $\zeta
_{n}\geq0,$ $\zeta_{n}^{\prime}\leq0$ and such that $\zeta_{n}\left(
\tau\right)  \rightarrow\chi_{\left[  \bar{\tau}_1, \bar \tau_2 \right]  } (\tau)$ as
$n\rightarrow\infty$. From (\ref{WeakFormSimp}), it then follows by means of a
standard computation that, for any function $\tilde{\psi}\in C^{1}\left(
\mathbb{R}_{*}\times\Delta^{d-1}\times\left[  \ln(2),\infty\right)  \right)  $
satisfying (\ref{TestLinEstim}),
\begin{align}
& \int_{\mathbb{R}_{*}\times \Delta^{d-1}}G\left(  \rho,\theta,\bar{\tau}%
_{2}\right)  \tilde{\psi}(\rho,\theta,\bar{\tau}_{2})d\Omega=\int_{\mathbb{R}_{*}\times \Delta^{d-1}}G\left(  \rho,\theta,\bar{\tau}_{1}\right)  \tilde{\psi
}(\rho,\theta,\bar{\tau}_{1})d\Omega\label{WeakSelfSimCut}\\
&  +\int_{\bar{\tau}_{1}}^{\bar{\tau}_{2}}d\tau\int_{\mathbb{R}_{*}\times \Delta^{d-1}}G\left(  \rho,\theta,\tau\right)  \left[  \partial_{\tau}%
\tilde{\psi}+\frac{\tilde{\psi}}{1-\gamma}-\frac{1}{1-\gamma}\rho
\partial_{\rho}\tilde{\psi}\right]  (\rho,\theta,\tau)d\Omega\nonumber\\
&  +\frac{1}{2}\int_{\bar{\tau}_{1}}^{\bar{\tau}_{2}}d\tau\int_{\mathbb{R}_{*}\times \Delta^{d-1}}\int_{\mathbb{R}_{*}\times \Delta^{d-1}}\tilde{K}(\rho
,\theta,r,\sigma,\tau)G\left(  \rho,\theta,\tau\right)  G\left(  r,\sigma
,\tau\right)  \Psi(\rho,r,\theta,\sigma,\tau)d\Omega d\tilde{\Omega}\nonumber
\end{align}
where $\tilde{K}$ and $\Psi$ are as in (\ref{Ktilde}) and (\ref{PsiTest}) respectively. 

We now use the test function $\tilde{\psi}(\rho,\theta,\tau)=\rho\Vert
\theta\Vert^{2}$ in (\ref{WeakSelfSimCut}). Notice that we have%
\[
\Psi(\rho,r,\theta,\sigma,\tau)=(\rho+r)\varphi(\rho+r,\frac{\rho}{\rho
+r}\theta,\frac{r}{\rho+r}\sigma)-\rho\varphi(\rho,\theta)-r\varphi
(r,\sigma)=-\frac{\rho r}{\rho+r}\Vert\theta-\sigma||^{2}.
\]

Using this identity in (\ref{WeakSelfSimCut}) we obtain%
\begin{align}
& \int_{\mathbb{R}_{*}\times \Delta^{d-1}}G\left(  \rho,\theta,\bar{\tau}%
_{2}\right)  \rho\Vert\theta\Vert^{2}d\Omega=
\int_{\mathbb{R}_{*}\times \Delta^{d-1}}G(\rho,\theta,\bar{\tau}_{1})\rho\Vert\theta\Vert^{2}%
d\Omega\nonumber\\
&  -\frac{1}{2}\int_{\bar{\tau}_{1}}^{\bar{\tau}_{2}}d\tau\int_{\mathbb{R}_{*}\times \Delta^{d-1}}\int_{\mathbb{R}_{*}\times \Delta^{d-1}}\tilde{K}(\rho
,\theta,r,\sigma,\tau)G\left(  \rho,\theta,\tau\right)  G\left(  r,\sigma
,\tau\right)  \frac{\rho r}{\rho+r}\Vert\theta-\sigma||^{2}d\Omega
d\tilde{\Omega} \label{IntLocEst}%
\end{align}
which combined with (\ref{eq:mass_conservation}) implies, taking $\bar{\tau
}_{1}=\ln(2),$ $\bar{\tau}_{2}=\bar{\tau},$ that%
\begin{equation}
\int_{\ln(2)}^{\bar{\tau}}d\tau\int_{\mathbb{R}_{*}\times \Delta^{d-1}}\int_{\mathbb{R}_{*}\times \Delta^{d-1}}\tilde{K}(\rho,\theta,r,\sigma,\tau)G\left(
\rho,\theta,\tau\right)  G\left(  r,\sigma,\tau\right)  \frac{\rho r}{\rho
+r}\Vert\theta-\sigma||^{2}d\Omega d\tilde{\Omega}\leq2 m_0 \label{IntLocFinite}%
\end{equation}
for all $\bar{\tau}\geq \ln(2).$

Using (\ref{eq:mass1}), (\ref{eq:sec_moment}) it follows that for any
$\delta_{0}>0$ small there exists $M>0$ (depending only on $\delta_{0}%
,\ m_{0}$ and $a$ in (\ref{eq:sec_moment})) such that%
\begin{equation}
\int_{(0,\frac{1}{M})\times\Delta^{d-1}}\rho G(\rho,\theta,\tau)d\Omega
+\int_{(M,\infty)\times\Delta^{d-1}}\rho G(\rho,\theta,\tau)d\Omega\leq
\delta_{0}m_{0}\ . \label{eq:delta}%
\end{equation}
Then%
\begin{equation}
\int_{\left[  \frac{1}{M},M\right]  \times\Delta^{d-1}}\rho G(\rho
,\theta,\tau)d\Omega\geq m_{0}\left(  1-\delta_{0}\right) . \label{eq:M}%
\end{equation}

The lower estimate in (\ref{KestimSimplex}) implies that there exists a
constant $\eta_{M}>0$ such that $\frac{\tilde{K}(\rho,\theta,r,\sigma,\tau
)}{\rho+r}\geq\eta_{M}$ for $\left(  \rho,r\right)  \in\left[  \frac{1}%
{M},M\right]  ^{2},\ \left(  \theta,\sigma\right)  \in\left(  \Delta
^{d-1}\right)  ^{2}$.
Using (\ref{IntLocFinite}) we then obtain
the estimate%
\begin{equation}
\int_{1}^{\bar{\tau}}d\tau\int_{\left[  \frac{1}{M},M\right] \times \Delta^{d-1}}\int_{\left[  \frac{1}{M},M\right]  \times \Delta^{d-1}}\rho
G\left(  \rho,\theta,\tau\right)  rG\left(  r,\sigma,\tau\right)  \Vert
\theta-\sigma||^{2}d\Omega d\tilde{\Omega}\leq\frac{2m_0 }{\eta_{M}} \label{LocMeasInt}%
\end{equation}
for all $\bar{\tau}\geq1 >\ln 2$.
For each $\tau\geq1$ and $M$ as above, we define a probability measure on
$\mathcal{M}_{+}\left(  \Delta^{d-1}\right)  $ by means of%
\begin{equation}
\lambda_{M}\left(  A,\tau\right)  =\frac{
\int_{\left[  \frac{1}%
{M},M\right] \times A  }\rho G\left(  \rho, \theta,\tau\right)
d\Omega}{\int_{\left[  \frac{1}%
{M},M\right] \times \Delta^{d-1}}\rho G\left(  \rho, \theta,\tau\right)
d\Omega}, \label{defLambda}%
\end{equation}
for each Borel set $A \subset\Delta^{d-1}$.
Then for each $\tau\geq1$ we have that $\lambda_{M}\left(  d\theta
,\tau\right)  $ is a probability measure and\ (\ref{eq:M}) and
(\ref{LocMeasInt}) imply, after taking the limit $\bar{\tau}\rightarrow\infty$%
\begin{equation}
\int_{1}^{\infty}d\tau\int_{\Delta^{d-1}}\int_{\Delta^{d-1}}\Vert\theta
-\sigma||^{2}\lambda_{M}\left(  d\theta,\tau\right)  \lambda_{M}\left(
d\sigma,\tau\right)  \leq\frac{2m_0}{\eta
_{M}\left(  m_{0}\right)  ^{2}\left(  1-\delta_{0}\right)  ^{2}}<\infty.
\label{dispEstimate}%
\end{equation}

It then follows that there exists a Borel set $\tilde{I}_{M}\subset\left(
1,\infty\right)  $ such that $\lim_{R\rightarrow\infty}\int_{\left(
R,\infty\right)  \cap\tilde{I}_{M}}d\tau=0$ and such that%
\[
\lim_{\tau\rightarrow\infty}\left[  \chi_{\left(  1,\infty\right)
\backslash\tilde{I}_{M}}\left(  \tau\right)  \int_{\Delta^{d-1}}\int
_{\Delta^{d-1}}\left\Vert \theta-\sigma\right\Vert ^{2}\lambda_{M}\left(
d\theta,\tau\right)  \lambda_{M}\left(  d\sigma,\tau\right)  \right]  =0
\]
where $\chi_{\left(  1,\infty\right)  \backslash\tilde{I}_{M}}$ is the
characteristic function of the set $\left(  1,\infty\right)  \backslash
\tilde{I}_{M}.$
We can then apply Lemma \ref{LemmSplit} to prove that there exists $\tau
_{0,M}$ sufficiently large such that, for any $\tau\in\left(  \tau
_{0,M},\infty\right)  \backslash\tilde{I}_{M}$ there exists a Borel set
$A_{M}\left(  \tau\right)  \subset\Delta^{d-1}$ such that the function defined
by means of%
\[
f_{M}\left(  \tau\right)  =\left\{
\begin{array}
[c]{c}%
\mathrm{diam}\left(  A_{M}\left(  \tau\right)  \right)  \ \ \text{for\ \ }%
\tau\in\left(  \tau_{0,M},\infty\right)  \backslash\tilde{I}_{M}\\
0\ \ \text{for\ \ }\tau\in\tilde{I}_{M}%
\end{array}
\right.
\]
satisfies $\lim_{\tau\rightarrow\infty}f_{M}\left(  \tau\right)  =0$ and, in
addition, the function%
\[
g_{M}\left(  \tau\right)  =\left\{
\begin{array}
[c]{c}%
\int_{A_{M}\left(  \tau\right)  }\lambda_{M}\left(  d\theta,\tau\right)
\ \ \text{for\ \ }\tau\in\left(  \tau_{0,M},\infty\right)  \backslash\tilde
{I}_{M}\\
1\ \ \text{for\ \ }\tau\in\tilde{I}_{M}%
\end{array}
\right.
\]
satisfies
\begin{equation}
\lim_{\tau\rightarrow\infty}g_{M}\left(  \tau\right)  =1. \label{locGlimit}%
\end{equation}

On the other hand, the conservation law (\ref{eq:mass_conservation}),
(\ref{eq:def_theta0}) combined with the fact that $\mathrm{diam}\left(
A_{M}\left(  \tau\right)  \right)  $ becomes arbitrarily small for any
$\tau\in\left(  \tau_{0,M},\infty\right)  \backslash\tilde{I}_{M}$ and
combined also with the estimate (\ref{eq:M}), implies that for $L>\tau_{0,M}$
large enough it holds 
\begin{equation}\label{eq:unionAM}
\bigcup_{\tau\in\left(  L,\infty\right)  \backslash
\tilde{I}_{M}}A_{M}\left(  \tau\right)  \subset B_{6\delta_{0}}\left(  \theta_{0}\right).
\end{equation}
Indeed, we have for each $\tau\in\left(  \tau_{0,M},\infty\right)  \backslash\tilde{I}_{M}$
\begin{align*}
\theta_{0}  &  =\frac{1}{m_{0}}\int_{{\mathbb{R}}_{*}\times\Delta^{d-1}}%
\rho\theta G(\rho,\theta,\tau)d\Omega\\
& =\frac{1}{m_{0}}\int_{\left[  \frac{1}{M},M\right]  \times\Delta^{d-1}%
}\rho\theta G(\rho,\theta,\tau) d\Omega+\frac{1}{m_{0}}\int_{\left(
{\mathbb{R}}_{*}\backslash\left[  \frac{1}{M},M\right]  \right)  \times
\Delta^{d-1}} \rho\theta G(\rho,\theta,\tau) d\Omega\\
&  =\frac{1}{m_{0}}\left[  \int_{\left[  \frac{1}{M},M\right] \times\Delta^{d-1}}\rho G\left(  \rho,\theta,\tau\right)  d\Omega\right]\int_{\Delta^{d-1}}\theta\lambda_{M}\left(  d\theta,\tau\right)\\
& \quad\quad +\frac{1}{m_{0}}\int_{\left(  {\mathbb{R}}_{*}\backslash\left[
\frac{1}{M},M\right]  \right)  \times\Delta^{d-1}}\rho\theta G(\rho,\theta,\tau
) d\Omega.
\end{align*}
Defining $m_{0,M}(\tau) := \int_{\left[  \frac{1}{M},M\right] \times\Delta^{d-1}}\rho G\left(  \rho,\theta,\tau\right)  d\Omega$ 
we obtain

\begin{align}\label{eq:theta0}
\theta_{0}  &   =\frac{m_{0,M}(\tau)}{m_{0}}
\int_{A_{M}\left(  \tau\right)  }\theta\lambda_{M}\left(  d\theta,\tau\right)
+ \frac{m_{0,M}(\tau)}{m_0} \int_{\Delta^{d-1}\backslash A_{M}(\tau)} \theta\lambda_M(d\theta,\tau)\nonumber \\
&+\frac{1}{m_{0}}\int_{\left(  {\mathbb{R}}_{*}\backslash\left[
\frac{1}{M},M\right]  \right)  \times\Delta^{d-1}}\rho\theta  G(\rho,\theta,\tau
) d\Omega.
\end{align}
We can write the third term on the right-hand side as $( 1 - \frac{m_{0,M}(\tau)}{m_0}),$ which can be estimated by $\delta_{0}$,  using (\ref{eq:delta}). Therefore, using (\ref{eq:M}),  \eqref{eq:theta0}, we obtain for  $\tau\in\left(  \tau_{0,M},\infty\right)  \backslash\tilde{I}_{M}$ sufficiently large, 
\begin{align*}
& \left\vert \int_{A_{M}\left(  \tau\right)  } \theta \lambda_{M}\left(  d\theta,\tau\right) -\theta_{0}\right\vert \leq \\
& \leq \left( 1 - \frac{m_{0,M}(\tau)}{m_0}\right)  \int_{A_{M}\left(  \tau\right)  } \theta \lambda_{M}\left(  d\theta,\tau\right)
 + (1- g_M(\tau))
+  \left(1-\frac{m_{0,M}(\tau)}{m_0}
\right)  \\
& \leq 2\delta_0 + (1- g_M(\tau)).
\end{align*}	
Thus using \eqref{locGlimit} we obtain, for $\tau$ sufficiently large, the estimate
\begin{align*}
\left\vert \int_{A_{M}\left(  \tau\right)  } \theta \lambda_{M}\left(  d\theta,\tau\right) -\theta_{0}\right\vert \leq   3\delta_0.
\end{align*}

Then, since $\mathrm{diam}\left(  A_{M}\left(  \tau\right)  \right)
\rightarrow0$ as $\tau\rightarrow\infty,$ and $\lambda_{M}\left(  d\theta
,\tau\right)  $ is a probability measure for each $\tau>\tau_{0,M},\ $ it
follows that for $\tau\geq L$ and $L$ sufficiently large and  $\tau\in\left(  \tau_{0,M},\infty\right)  \backslash\tilde{I}_{M}$ we have 
\[
\left\vert \int_{A_{M}\left(  \tau\right)  }\theta\lambda_{M}\left(  d\theta
,\tau\right)   -\bar{\theta}\left(  \tau\right)  \right\vert
\leq\delta_{0}%
\]
for some $\bar{\theta}\left(  \tau\right)  \in\Delta^{d-1}.$ 
Then, combining the last two previous inequalities, we obtain $\left\vert
\bar{\theta}\left(  \tau\right)  -\theta_{0}\right\vert \leq 4\delta_{0},$ for $\tau\in\left(  \tau_{0,M},\infty\right)  \backslash\tilde{I}_{M}$ sufficiently large. We
have also $A_{M}\left(  \tau\right)  \subset B_{\mathrm{diam}\left(
A_{M}\left(  \tau\right)  \right)   +\delta_0}\left(  \bar{\theta}\left(  \tau\right)
\right)  \subset B_{2\delta_0}\left(  \bar{\theta}\left(  \tau\right)
\right),$ if $\tau \geq L$ and $L$ sufficiently large. In addition we obtain $A_{M}\left(
\tau\right)  \subset B_{6\delta_{0}}\left(  \theta_{0}\right)  $ for any
$\tau\in\left(  L,\infty\right)  \backslash\tilde{I}_{M},$ hence the claim \eqref{eq:unionAM} follows.

In order to conclude the proof we consider a sequence of values $\delta
_{0}=\frac{1}{n},\ n\in\mathbb{N}$, the corresponding sequence $M_{n}%
\rightarrow\infty,$ the sets $\tilde{I}_{M_{n}}$ and a sequence of increasing
values $\tau_{0,M_{n}}$ such that $\left\vert \left(  \tau_{0,M_{n}}%
,\infty\right)  \cap\tilde{I}_{M_{n}}\right\vert \leq2^{-n}$ for
$n\in\mathbb{N}$ with $\lim_{n\rightarrow\infty}\tau_{0,M_{n}}=\infty$ and
also $\mathrm{diam}\left(  A_{M_{n}}\left(  \tau\right)  \right)  \leq2^{-n}$
for $\tau\in\left(  \tau_{0,M_{n}},\tau_{0,M_{n+1}}\right]  \backslash
\tilde{I}_{M_{n}}$ and $\int_{A_{M_{n}}\left(  \tau\right)  }\lambda_{M_{n}%
}\left(  d\theta,\tau\right)  \geq1-2^{-n}$ for $\tau\in\left(  \tau_{0,M_{n}%
},\tau_{0,M_{n+1}}\right]  \backslash\tilde{I}_{M_{n}}.$ Then, if we define
$\tilde{I}=\bigcup_{n\in\mathbb{N}}\tilde{I}_{M_{n}}$ it follows that
$\lim_{L\rightarrow\infty}\left(  L,\infty\right)  \cap\tilde{I}=0.$

We define the sets%
\[
A\left(  \tau\right)  =A_{M_{n}}\left(  \tau\right)  \ \ \text{for\ \ }\tau
\in\left(  \tau_{0,M_{n}},\tau_{0,M_{n+1}}\right]  \backslash\tilde{I}.
\]
We have also that the function defined as%
\[
f\left(  \tau\right)  =\left\{
\begin{array}
[c]{c}%
\mathrm{diam}\left(  A\left(  \tau\right)  \right)  \ \ \text{for\ \ }\tau
\in\left(  1,\infty\right)  \backslash\tilde{I}\\
0\ \ \text{for\ \ }\tau\in\tilde{I}%
\end{array}
\right.
\]
satisfies $\lim_{\tau\rightarrow\infty}f\left(  \tau\right)  =0.$ Moreover,
the function $g$ defined as%
\[
g\left(  \tau\right)  =\left\{
\begin{array}
[c]{c}%
\int_{A\left(  \tau\right)  }\lambda_{M_{n}}\left(  d\theta,\tau\right)
\ \ \text{for\ \ }\tau\in\left(  \tau_{0,M_{n}},\tau_{0,M_{n+1}}\right]
\backslash\tilde{I}\\
1\ \ \text{for\ \ }\tau\in\tilde{I}%
\end{array}
\right.
\]
satisfies $\lim_{\tau\rightarrow\infty}g\left(  \tau\right)  =1.$ It then
follows that $\bigcap_{\left\{  L\geq1\right\}  }\bigcup_{\tau\in\left(
L,\infty\right)  \backslash\tilde{I}_{M}}A\left(  \tau\right)  =\left\{
\theta_{0}\right\}  .$ On the other hand, the definition of the function
$g\left(  \cdot\right)  $ implies, as well as the definition of the measures
$\lambda_{M_{n}}$ (c.f. (\ref{defLambda})) implies (\ref{locEstimate}) for
some function $\delta\left(  t\right)  $ such that $\lim_{t\rightarrow\infty
}\delta\left(  t\right)  =0.$ Using the fact that $d\tau=\frac{dt}{t+1}$ it
follows that if we define the set $I\subset\left(  0,\infty\right)  $ as the
image of $\tilde{I}$ by means of the mapping $\tau\rightarrow\left(  e^{\tau
}-1\right)  $, it then follows from $\int_{L}^{\infty}d\tau\rightarrow0$ as
$L\rightarrow\infty$ that $\lim_{T\rightarrow\infty}\frac{\left\vert
I\cap\left[  L,2L\right]  \right\vert }{L}=0.$ This concludes the proof of
Theorem \ref{mainRough}.
\end{proofof}


We now prove Theorem \ref{mainRoughImproved}.


\begin{proofof}
[Proof of Theorem \ref{mainRoughImproved}]The proof is similar to the one of
Theorem \ref{mainRough}. The main difference is that under the assumptions in
Theorem \ref{mainRoughImproved} we can prove that the measure $\lambda
_{M}\left(  d\theta,\tau\right)  $ changes continuously in the weak$-\ast$
topology as $\tau$ varies. More precisely, it turns out that if $0\leq
\gamma<1$ and $0\leq\gamma+p<1$ the following estimate holds%
\begin{equation}
\int_{{\mathbb{R}}_{*}\times\Delta^{d-1}}G(\rho,\theta,\tau)\rho^{\sigma
}d\Omega\leq C_{0}\ \ ,\ \ \tau\geq1 \label{MomentSigma}%
\end{equation}
where $\sigma=\gamma$ if $p>0$ and $\sigma=\gamma+\delta$ with $\delta>0$
arbitrarily small if $p\leq0.$ The constant $C_{0}$ on the right-hand side of
(\ref{MomentSigma}) is independent of $\tau.$

 The estimate (\ref{MomentSigma}) has been proved in \cite{EM06} in the case of
one-component coagulation equations (i.e. $d=1\,)$ under the assumption
$f_{0} \in L_{loc}^1\left({\mathbb{R}}_{\ast}^{d}\right)$. The proofs can be adapted to the case in which $f_0 \in \mathcal{M}_+(\R_*^d)$. A sketch of the ideas required to prove this estimate in the multicomponent coagulation case $d>1$
are collected in Section \ref{MomEst} (cf.\ Proposition \ref{PropMomEst}).

On the other hand we also have the conservation of mass identity (cf.
(\ref{eq:mass1}))%
\begin{equation}
\int_{{\mathbb{R}}_{*}\times\Delta^{d-1}}G(\rho,\theta,\tau)\rho d\Omega
=m_{0}\ \ ,\ \ \tau\geq1. \label{GMass}%
\end{equation}

Using (\ref{MomentSigma}) and (\ref{GMass}), as well as the fact that for the
range of parameters under consideration we have $\gamma<1-p$ it follows that
\begin{equation}
\int_{\left\{  \rho\geq1\right\}  \times\Delta^{d-1}}G(\rho,\theta,\tau
)\rho^{\gamma+p}d\Omega+\int_{\left\{  \rho\leq1\right\}  \times\Delta^{d-1}%
}G(\rho,\theta,\tau)\rho^{1-p}d\Omega\leq C_{1}\ \ ,\ \ \tau\geq1
\label{Gmoment}%
\end{equation}
where $C_{1}$ depends on $C_{0}$ and $m_{0}$ but it is independent of $\tau.$
Using (\ref{WeakSelfSimCut}) we obtain that, for any smooth test function
$\tilde{\psi}$ satisfying (\ref{TestLinEstim}) and any $\bar\tau_{2} \geq
\bar\tau_{1} >0$ we have%
\begin{align}
& \int_{\mathbb{R}_{*}\times \Delta^{d-1}}G\left(  \rho,\theta,\bar{\tau}%
_{2}\right)  \tilde{\psi}(\rho,\theta,\bar{\tau}_{2})d\Omega-\int_{\mathbb{R}_{*}\times \Delta^{d-1}}G\left(  \rho,\theta,\bar{\tau}_{1}\right)  \tilde{\psi
}(\rho,\theta,\bar{\tau}_{1})d\Omega\nonumber\\
&  =\int_{\bar{\tau}_{1}}^{\bar{\tau}_{2}}d\tau\int_{ \mathbb{R}_{*}\times\Delta^{d-1}}G\left(  \rho,\theta,\tau\right)  \left[  \partial_{\tau}%
\tilde{\psi}+\frac{\tilde{\psi}}{1-\gamma}-\frac{1}{1-\gamma}\rho
\partial_{\rho}\tilde{\psi}\right]  (\rho,\theta,\tau)d\Omega\nonumber\\
&  +\frac{1}{2}\int_{\bar{\tau}_{1}}^{\bar{\tau}_{2}}d\tau\int_{ \mathbb{R}_{*}\times\Delta^{d-1}}\int_{ \mathbb{R}_{*}\times\Delta^{d-1}}\tilde{K}(\rho
,\theta,r,\sigma,\tau)G\left(  \rho,\theta,\tau\right)  G\left(  r,\sigma
,\tau\right)  \Psi(\rho,r,\theta,\sigma,\tau)d\Omega d\tilde{\Omega}.
\label{DiffEstim}%
\end{align}

We now argue as in the proof of Theorem \ref{mainRough}. Given $\delta_{0}>0$
arbitrarily small we select $M>0$ sufficiently large such that (\ref{eq:delta}%
), (\ref{eq:M}) hold. We then define $\lambda_{M}\left(  d\theta,\tau\right)
$ by means of\ (\ref{defLambda}). Then (\ref{dispEstimate}) holds.

Suppose that $\tilde{\psi}(\cdot,\cdot, \tau) \in W^{1,\infty}\left(  \R_*
\times\Delta^{d-1} \right)  $ such that $\operatorname*{supp}\left(
\tilde{\psi}(\cdot,\cdot, \tau)\right)  \subset\R_*  \times\Delta^{d-1}$ and
$\operatorname*{supp}\left(  \tilde{\psi}\right) $  is compact. Then, using
(\ref{KestimSimplex}), (\ref{GMass}) and (\ref{Gmoment}), we obtain from
\eqref{DiffEstim} the following estimate%
\begin{align*}
\left\vert\int_{\mathbb{R}_{*}\times \Delta^{d-1}}G\left(  \rho,\theta,\bar{\tau
}_{2}\right)  \tilde{\psi}(\rho,\theta,\bar{\tau}_{2})d\Omega-\int_{ \mathbb{R}_{*}\times\Delta^{d-1}}G\left(  \rho,\theta,\bar{\tau}_{1}\right)
\tilde{\psi}(\rho,\theta,\bar{\tau}_{1})d\Omega\right\vert \leq\\
C_{2}\max_{\tau \geq 1 }\left\{\left\Vert \tilde{\psi}(\cdot,\cdot, \tau)\right\Vert _{W^{1,\infty}\left( \R_* \times\Delta^{d-1}\right)  }\right\} \left\vert \bar{\tau}_{2}%
-\bar{\tau}_{1}\right\vert
\end{align*}
for $\bar{\tau}_{1},\bar{\tau}_{2}\geq1$, where $C_{2}$ depends on $m_{0}$ and
$C_{0}$. It then follows that the mapping  $\tau\in\left[  1,\infty\right)
\to G\left( \rho,\theta,\tau\right) \in\mathcal{M}_{+}\left( \mathbb{R}%
_{+}\times\Delta^{d-1}\right)  $ is continuous in the weak$-\ast$ topology of
$\mathcal{M}_{+}\left(  \mathbb{R}_{+}\times\Delta^{d-1}\right) .$ In
particular, this implies that the mapping $\tau\rightarrow\int_{\left[
\frac{1}{M},M\right]  }\int_{\Delta^{d-1}}\rho G\left(  \rho,\theta
,\tau\right)  d\Omega$ is continuous in $\tau\geq1$ and that the mapping
$\tau\in\left[  1,\infty\right)  \to\lambda_{M}\left(  d\theta,\tau\right)
\in\mathcal{M}_{+}\left(  \Delta^{d-1}\right)  $ with $\lambda_{M}$ as in
(\ref{defLambda}) is also continuous in the weak$-\ast$ topology of
$\mathcal{M}_{+}\left(  \Delta^{d-1}\right)  .$ Moreover, the mapping from
$\mathcal{M}_{+}\left(  \Delta^{d-1}\right)  $ to $\left[  0,\infty\right)  $
defined by means of $\lambda_{M}\rightarrow\int_{\Delta^{d-1}}\int
_{\Delta^{d-1}}\Vert\theta-\sigma||^{2}\lambda_{M}\left(  d\theta\right)
\lambda_{M}\left(  d\sigma\right)  $ is continuous if the topology of
$\mathcal{M}_{+}\left(  \Delta^{d-1}\right)  $ is given by the weak$-\ast$
topology. This follows from the fact that the tensor product is a continuous
mapping in the weak$-\ast$ topology. We now claim that (\ref{dispEstimate})
implies that
\begin{equation}
\lim_{\tau\rightarrow\infty}\int_{\Delta^{d-1}}\int_{\Delta^{d-1}}\Vert
\theta-\sigma||^{2}\lambda_{M}\left(  d\theta,\tau\right)  \lambda_{M}\left(
d\sigma,\tau\right)  =0. \label{limLambda}%
\end{equation}
Indeed, suppose that $\lim\sup_{\tau\rightarrow\infty}\int_{\Delta^{d-1}}%
\int_{\Delta^{d-1}}\Vert\theta-\sigma||^{2}\lambda_{M}\left(  d\theta
,\tau\right)  \lambda_{M}\left(  d\sigma,\tau\right)  >0.$ Then, there exist
an increasing sequence $\left\{  \tau_{n}\right\}  _{n\in\mathbb{N}}$ with
$\lim_{n\rightarrow\infty}\tau_{n}=\infty$ and $\eta>0$ such that%
\[
\int_{\Delta^{d-1}}\int_{\Delta^{d-1}}\left\Vert \theta-\sigma\right\Vert
^{2}\lambda_{M}\left(  d\theta,\tau_{n}\right)  \lambda_{M}\left(
d\sigma,\tau_{n}\right)  >\eta
\]
for $n$ large enough. We can assume without loss of generality that
$\tau_{n+1}-\tau_{n}\geq1.$ Then, the uniform continuity of the mapping
$\tau\rightarrow\int_{\Delta^{d-1}}\int_{\Delta^{d-1}}\left\Vert \theta
-\sigma\right\Vert ^{2} \lambda_{M}\left(  d\theta,\tau\right)  \lambda
_{M}\left(  d\sigma,\tau\right)  $ implies that there exists $\varepsilon
_{0}>0$ small such that
\[
\int_{\Delta^{d-1}}\int_{\Delta^{d-1}}\left\Vert \theta-\sigma\right\Vert ^{2}
\lambda_{M}\left(  d\theta,\tau\right)  \lambda_{M}\left(  d\sigma
,\tau\right)  >\frac{\eta}{2}\text{ for }\tau\in\left(  \tau_{n}%
-\varepsilon_{0},\tau_{n}+\varepsilon_{0}\right) .
\]

However, this contradicts (\ref{dispEstimate}) and implies
(\ref{limLambda}). We can then apply Lemma \ref{LemmSplit} to show that there
exists a family of Borel sets $\left\{  A_{M}\left(  \tau\right)  \right\}
_{\tau\geq1}$ with $A_{M}\left(  \tau\right)  \subset\Delta^{d-1}$ and
$\lim_{\tau\rightarrow\infty}\mathrm{diam}\left(  A_{M}\left(  \tau\right)
\right)  =0$ such that%
\begin{equation}
\lim_{\tau\rightarrow\infty}\int_{A_{M}\left(  \tau\right)  }\lambda
_{M}\left(  d\theta,\tau\right)  =1. \label{locGlimRest}%
\end{equation}

Notice that this result is similar to the claim (\ref{locGlimit}) in the proof
of Theorem \ref{mainRough}, with the only difference that in the case of
(\ref{locGlimit}) there is a set of\ ``exceptional" times with small
measure for which $\int_{A_{M}\left(  \tau\right)  }\lambda_{M}\left(
d\theta,\tau\right)  $ might not be close to $1.$ We can now take $\delta
_{0}=\frac{1}{n},$ $n\in\mathbb{N}$, select the corresponding values of the sequence 
$M_{n}\to \infty$ as it was made in the proof of Theorem \ref{mainRough} and argue
exactly as it was made there in order to prove (\ref{localRestricted}). Hence 
the result follows.
\end{proofof}


\subsection{Complete localization along a ray in self-similar solutions}


We now prove Theorem \ref{TheorSelfSimLoc}.


\begin{proofof}
[Proof of Theorem \ref{TheorSelfSimLoc}]Using the change of variables
(\ref{rhoTheta}), (\ref{PolDist}) with $F$ and $G$ independent of time, we can
rewrite (\ref{WeakSelfSim}) as%
\begin{align}
&  \frac{1}{2}\int_{ \mathbb{R}_{*}\times\Delta^{d-1}}\int_{ \mathbb{R}_{*}\times\Delta^{d-1}}\hat{K}(\rho,\theta,r,\sigma)G\left(  \rho,\theta\right)
G\left(  r,\sigma\right)  \Psi(\rho,r,\theta,\sigma)d\Omega d\tilde{\Omega
}\label{WeakSelfSimRho}\\
&  +\frac{1}{1-\gamma}\int_{ \mathbb{R}_{*}\times\Delta^{d-1}}G\left(  \rho
,\theta\right)  \left[  \tilde{\psi}-\rho\partial_{\rho}\tilde{\psi}\right]
(\rho,\theta)d\Omega  =0\nonumber
\end{align}
where due to the homogeneity of the kernel $K$ we have%
\[
\hat{K}(\rho,\theta,r,\sigma)=K(\rho\theta,r\sigma)\ ,\ \ \rho,r\in
\mathbb{R}_{+}\ ,\ \ \theta,\sigma\in\Delta^{d-1}%
\]
and%
\begin{equation}
\Psi(\rho,r,\theta,\sigma)=\tilde{\psi}\left(  \rho+r,\frac{\rho}{\rho
+r}\theta+\frac{r}{\rho+r}\sigma\right)  -\tilde{\psi}(\rho,\theta
)-\tilde{\psi}(r,\sigma).
\end{equation}
Notice that $\tilde{\psi}(\rho,\theta)=\psi\left(  \xi\right)  \ $with
$\xi=\rho\theta.$

Arguing as in the proof of Theorem \ref{mainRough} we can prove that
(\ref{WeakSelfSimRho}) holds for any test function $\tilde{\psi}(\rho,\theta)$
satisfying (cf.\ (\ref{TestLinEstim}))%
\[
\left\vert \tilde{\psi}(\rho,\theta)\right\vert +\rho\left\vert \frac
{\partial\tilde{\psi}}{\partial\rho}(\rho,\theta)\right\vert +\left\vert
\nabla_{\theta}\tilde{\psi}(\rho,\theta)\right\vert \leq C\rho\ \ ,\ \ \rho
\in\mathbb{R}_{+},\ \theta\in\Delta^{d-1}.
\]

We can then choose the test function $\tilde{\psi}(\rho,\theta)=\rho\left\Vert
\theta\right\Vert ^{2}$ in (\ref{WeakSelfSimRho}). Then%
\begin{equation}
\int_{ \mathbb{R}_{*}\times\Delta^{d-1}}\int_{ \mathbb{R}_{*}\times\Delta^{d-1}}
\hat{K}(\rho,\theta,r,\sigma)G\left(  \rho,\theta\right)  G\left(
r,\sigma\right)  \frac{\rho r}{\rho+r}\Vert\theta-\sigma||^{2}d\Omega
d\tilde{\Omega}=0. \label{GIdentZero}%
\end{equation}

This identity implies that $G$ has the form%
\begin{equation}
G\left(  \rho,\theta\right)  =G_{0}\left(  \rho\right)  \delta\left(
\theta-\theta_{0}\right) \label{GSelfSimDirac}%
\end{equation}
where $G_{0} \in\mathcal{M}_{+}({\mathbb{R}}_{*} )$.

This can be seen defining for each $M$ sufficiently large the probability
measures on $\Delta^{d-1}$ by
\begin{equation}
\lambda_{M}( A) =\frac{\int_{\left[  \frac{1}{M},M\right] \times A} \rho G\left(
\rho,\theta\right)  d\Omega}{\int_{\left[  \frac{1}{M},M\right]  \times\Delta^{d-1}}\rho G\left(  \rho,\theta\right)  d\Omega} \label{LambMDef}%
\end{equation}
for each Borel set $A \subset\mathcal{M}_{+}(\Delta^{d-1})$.

These probability measures are well defined for $M$ large enough since $G$ is
not identically zero. Then (\ref{GIdentZero}) implies%
\[
\int_{\Delta^{d-1}}\int_{\Delta^{d-1}}\lambda_{M}\left(  d\theta\right)
\lambda_{M}\left(  d\sigma\right)  \Vert\theta-\sigma||^{2}=0.
\]
We can then apply Lemma \ref{LemmSplitDelta} with $\varepsilon$ and $\delta$
arbitrarily small to prove that%
\begin{equation}
\lambda_{M}
= \delta_{\theta_{M}}
\label{LambDelta}%
\end{equation}
with $\theta_{M}\in\Delta^{d-1}.$ We now use that, since $\bigcup_{M>1}\left[
\frac{1}{M},M\right]  \times\Delta^{d-1}=\mathbb{R}_{\ast}^{d}$
\[
\int_{\left[  \frac{1}{M},M\right]  \times\Delta^{d-1}}G\left(  \rho
,\theta\right)  \rho d\Omega\rightarrow m_{0}\text{ \ as }M\rightarrow\infty.
\]

Hence, combining the previous limit with \eqref{LambDelta} and using the
change of variables (\ref{rhoTheta}), (\ref{PolDist}), the definition of
$\theta_{0}$ in \eqref{DefThetaZero} and Lebesgue dominated convergence
Theorem, we obtain
\[
\lim_{M\rightarrow\infty}\int_{\Delta^{d-1}}\theta\lambda_{M}\left(
d\theta\right)  =\frac{1}{m_{0}}\int_{ \mathbb{R}_{*}\times \Delta^{d-1}}G\left(
\rho,\theta\right)  \theta\rho d\Omega=\theta_{0}.
\]

Therefore, from (\ref{LambDelta}) we obtain that $\lim_{M\rightarrow\infty
}\theta_{M}=\theta_{0}.$ Notice that (\ref{LambMDef}) and (\ref{LambDelta})
imply that there exists a measure $G_{0,M}\in\mathcal{M}_{+}\left(
\R_*\right)  $ such that%
\[
G\left(  \rho,\theta\right)  \chi_{\left[  \frac{1}{M},M\right]  }\left(
\rho\right)  =G_{0,M}\left(  \rho\right)  \delta\left(  \theta-\theta
_{M}\right) .
\]

This implies that the measure $G\left(  \rho,\theta\right)  $ is supported along
the line $\left\{  \theta=\theta_{M}\right\}  .$ Therefore $\theta_{M}$ is
independent of $M$ and we have $\theta_{M}=\theta_{0}.$ This gives
(\ref{GSelfSimDirac}).

Plugging (\ref{GSelfSimDirac}) into (\ref{WeakSelfSimRho}) we obtain that
$G_{0}$ satisfies
\begin{align*}
&  \frac{1}{2}\int_{ \mathbb{R}_{*}}\int_{ \mathbb{R}_{*}}\hat{K}(\rho,\theta
_{0},r,\theta_{0})G_{0}\left(  \rho\right)  G_{0}\left(  r\right)  \Psi
(\rho,r,\theta_{0},\theta_{0}) r^{d-1}\rho^{d-1} d\rho dr\\
&  +\frac{1}{1-\gamma}\int_{ \mathbb{R}_{*}}G_{0}\left(  \rho\right)  \left[
\tilde{\psi}-\rho\partial_{\rho}\tilde{\psi}\right]  (\rho,\theta_{0}
)\rho^{d-1}d\rho  =0
\end{align*}
for each $\tilde{\psi}\in C_{c}^{1}\left(  \left(  0,\infty\right)
\times\Delta^{d-1}\right)  .$ In particular, defining $F_0:=\rho^{d-1} G_{0}$ and $K_{\theta_0}(\rho,r) = \hat{K}(\rho,\theta_{0},r,\theta_{0})$, we obtain
\begin{align*}
&  \frac{1}{2}\int_{ \mathbb{R}_{*}}\int_{ \mathbb{R}_{*}} K_{\theta_0}(\rho,r)F_{0}\left(  \rho\right)  F_{0}\left(  r\right)  \left[
\tilde{\psi}\left(  \rho+r\right)  -\tilde{\psi}(\rho)-\tilde{\psi}(r)\right]
 d\rho dr\\
&  +\frac{1}{1-\gamma}\int_{ \mathbb{R}_{*}}F_{0}\left(  \rho\right)  \left[
\tilde{\psi}-\rho\partial_{\rho}\tilde{\psi}\right]  (\rho) d\rho  =0
\end{align*}
for any $\tilde{\psi}\in C_{c}^{1}\left( 0,\infty\right)  ,$ i.e. $F_0$ is a self-similar profile for the
one-dimensional coagulation equation with coagulation kernel $K_{\theta_0}(\rho,r).$ This concludes the proof of Theorem
\ref{TheorSelfSimLoc}.

\end{proofof}


\section{Global existence and self-similar solutions for the multicomponent
problem\label{GlobEx}}


In this Section we show that the assumptions of Theorems 
\ref{mainRoughImproved} and \ref{thm:concentDisc} are satisfied for some ranges of exponents
$\gamma,\ p$ as well as for some choices of initial values $f_{0}%
\in\mathcal{M}_{+}({\mathbb{R}}_{\ast}^{d})$. The global existence of weak
solutions in the sense of  Definition \ref{def:time-dep_sol} which in
addition satisfy the conservation of mass condition (\ref{MassCons}) has been
proved (cf.\ \cite{book,EM06, FL06,Norris}) in the case of one-component
coagulation equations (i.e. $d=1$) for product kernels of the form
$x^{\gamma+\lambda} y^{-\lambda} +x^{\gamma+\lambda} y^{-\lambda}$ with $-\lambda \leq\gamma+\lambda< 1$ and
$\gamma<1$. These results can be easily extended to the class of kernels
considered here satisfying \eqref{eq:condK_cont}, (\ref{eq:condK_sym2}), (\ref{eq:Phi_kernel}),
(\ref{KernIneq}),  with $d \geq1$.


On the other hand, the existence of self-similar profiles has been proved for
one-component coagulation kernels which in addition to the previously stated
conditions (cf.\ (\ref{eq:condK_sym2}), (\ref{eq:Phi_kernel}), (\ref{KernIneq}%
), \eqref{eq:condK_cont}) satisfy also (\ref{eq:condK_sym1}) and the
homogeneity condition (\ref{eq:condK_homogeneity}). These results can be used
to prove the existence of self-similar solutions with the form
(\ref{SelfSimLoc}) in the multicomponent case (i.e. $d>1$). We will explain in
this Section how this can be achieved.


We first notice that the following global existence result holds.


\begin{theorem}
\label{thm:existence_td} Suppose that $K$ is as in (\ref{eq:condK_cont}) and
satisfies  the homogeneity property (\ref{eq:condK_homogeneity}) and the upper and lower bounds (\ref{eq:condK_sym1}), \eqref{eq:condK_sym2}
 with $\gamma,\ p\in\mathbb{R}$
such that 
 $ \gamma,\gamma+p \in [0,1)$. Suppose that  $f_0 \in \mathcal{M}_+(\R_*^d)$ satisfies
\begin{equation}
\int_{{\mathbb{R}}_{*}^{d}}(|x|+|x|^{1+r})f_{0}(dx)<\infty\label{eq:ini_cond}%
\end{equation}
for some $r>0$. Then, there exists a weak solution $f\in C([0,\infty
),\mathcal{M}_{+}({{\mathbb{R}}_{\ast}^{d}}))$ to \eqref{B4} in the
sense of Definition \ref{def:time-dep_sol} with $f(0,\cdot)=f_{0}$.
 Moreover, this solution $f$ has the following property. For any $k$ satisfying $k \in(\gamma,1+r]$ if $p\leq0$ or $k
\in(-\infty, 1+r]$ if $p>0$, there is a constant $c>0$ that may depend on the
initial data such that
\begin{equation}
\int_{{\mathbb{R}}_{*}^{d} }|x|^{k}f(dx,t) \leq c t^{\frac{k-1}{1-\gamma}%
},\ \quad t\geq1 . \label{MomentEstThm4.1}
\end{equation}

\end{theorem}

The existence of a solution under the conditions of Theorem
\ref{thm:existence_td} can be obtained by adapting to the multicomponent setting
the results proved in \cite{EM06} (see also \cite{book}).  Those results
have been obtained for locally integrable initial data $f_{0}\in
{L}^{1}_{\text{loc}}({\mathbb{R}}_{\ast}^{d})$, but they can be adapted to the case of more general initial data $f_{0}\in\mathcal{M}_{+}({{\mathbb{R}}_{\ast}^{d}}).$  The moment estimate
\eqref{MomentEstThm4.1} has been derived in \cite{EM06} for $d=1$ (cf.\ Theorem $2.4$ and Lemma $3.1$ in \cite{EM06}). We will present in Section \ref{MomEst} the ideas that allow to generalize 
the proof to multicomponent coagulation
equations.

On the other hand, the existence of self-similar profiles is well known for a
large class of homogeneous kernels $K$ in one component (i.e. $d=1$)
coagulation systems. Using the results obtained for the one component system
in \cite{EM05, EM06, FL04} we can immediately prove the existence of
self-similar profiles for the multicomponent system in terms of the solution
to the one-component equation, having the form (\ref{SelfSimLoc}). Moreover,
Theorem \ref{TheorSelfSimLoc} guarantees that all  self-similar profiles of
(\ref{B4}) in the sense of Definition \ref{def:self_sim_sol} have the
form (\ref{SelfSimLoc}). We have the following result.

\begin{theorem}
\label{thm:existence_self_sim}Suppose that $K:({{\mathbb{R}}_{\ast}^{d}})^{2}%
\rightarrow{\mathbb{R}}_{+}$ is as in (\ref{eq:condK_cont}) and satisfies the homogeneity property (\ref{eq:condK_homogeneity}) and the upper
and lower bounds (\ref{eq:condK_sym1}), \eqref{eq:condK_sym2} with
  $ \gamma,\gamma+p \in [0,1)$. Let $m\in{{\mathbb{R}}_{\ast}^{d}}$
with the form $m=\left(  m_{k}\right)  _{k=1}^{d}$ with $m_{k}\geq0$ for any
$k=1,2,...,d$ and satisfying $\left\vert m\right\vert >0.$ Let $\theta_0
=\frac{m}{\left\vert m\right\vert }\in\Delta^{d-1}.$ There exists at least one
measure $F_{0}\in\mathcal{M}_{+}({{\mathbb{R}}_{\ast}})$ such that the measure
$F\in\mathcal{M}_{+}({{\mathbb{R}}_{\ast}^{d}})$ defined in (\ref{SelfSimLoc})
is a self-similar profile to \eqref{B4} in the sense of Definition
\ref{def:self_sim_sol}. Moreover, we have
\begin{equation}
\int_{{{\mathbb{R}}_{\ast}^{d}}}\xi_{k}F(d\xi)=m_{k}\ ,\ \ k=1,2,...,d.
\label{eq:mass_self_sim}%
\end{equation}

\end{theorem}

\begin{proof}
Suppose that $F$ has the form (\ref{SelfSimLoc}) with $F_{0}\in\mathcal{M}
_{+}({{\mathbb{R}}_{\ast}})$ a solution to the one-dimensional problem. We
then have, using the variables $\left(  \rho,\theta\right) $, that
\[
F(d\xi)=F_{0}\left(  \rho\right) \delta\left(
\theta-\theta_{0}\right)  d\rho d\nu\left(  \theta\right)  .
\]
Hence, (\ref{eq:bound_gamma_moment_selfsim}) holds if and only if%
\begin{equation}
\int_{\left(  1,\infty\right)  }\rho^{\gamma+p}F_{0}\left(  \rho\right)
d\rho+\int_{\left(  0,1\right]  }\rho^{-p+1}F_{0}\left(  \rho\right)
d\rho+\int_{{{\mathbb{R}}_{\ast}}}\rho F_{0}\left(  \rho\right)
d\rho<\infty\label{EstSelfSimF_0}%
\end{equation}
and (\ref{WeakSelfSim}) holds if and only if the following identity is
satisfied%
\begin{align}
&  \int_{{{\mathbb{R}}_{\ast}}}\int_{{{\mathbb{R}}_{\ast}}}K\left(  \rho
\theta_{0},r\theta_{0}\right)  \rho\left[  \varphi(\rho+r,\theta_{0})-\varphi
(\rho,\theta_{0})\right]  F_{0}\left(  \rho\right)  F_{0}\left(  r\right)
d\rho dr\label{WeakSelfSimF_0}\\
&  -\frac{1}{1-\gamma}\int_{{{\mathbb{R}}_{\ast}}}F_{0}\left(  \rho\right)
\frac{\partial\varphi}{\partial\rho}(\rho,\theta_{0})\rho^{2}d\rho\nonumber\\
&  =0\nonumber
\end{align}
with $\psi(\xi) = \rho \varphi
(\rho,\theta),\ \xi=\rho\theta$ and $\varphi\in C_{c}%
^{1}({{\mathbb{R}}_{\ast}^{d}})$ is an arbitrary test function. Notice that since
(\ref{EstSelfSimF_0}) holds, then using the definition of $\theta_0$, (\ref{eq:mass_self_sim}) is 
automatically satisfied.

We define $K_{\theta_0}\left(
\rho,r\right)  =K\left(  \rho\theta_{0},r\theta_{0}\right)$. 
Notice that the kernel $K_{\theta_0}$ is homogeneous with homogeneity $\gamma$ and continuous. Due to (\ref{eq:condK_sym1}), \eqref{eq:condK_sym2} we
have that $K_{\theta_0}$ satisfies
\begin{equation}
c_{1}(\rho+r)^{\gamma}\Phi_{p}\left( \frac{\rho}{\rho+r}\right) \leq K_{\theta_0}\left(  \rho,r\right)  \leq c_{2}(\rho+r)^{\gamma}\Phi_{p}\left( \frac
{\rho}{\rho+r}\right) \ \ ,\ \ \rho,r\in{{\mathbb{R}}_{\ast}}.
\label{KtildaEst}%
\end{equation}
The existence of measures $F_{0}\in\mathcal{M}_{+}({{\mathbb{R}}_{\ast}})$
satisfying (\ref{WeakSelfSimF_0}), (\ref{EstSelfSimF_0}) for kernels
satisfying (\ref{KtildaEst}) with $\gamma,\ p$ satisfying  $\gamma,\gamma+p \in [0,1)$  is ensured by the results in \cite{EM05, EM06,
FL04}. Then the result follows.
\end{proof}

\section{Moment estimates\label{MomEst}}

A crucial step in the proof of the existence of self-similar solutions for
one-dimensional coagulation equations is the derivation of some estimates for
the moments of $f$ which guarantee that the mass of the monomers of the
solutions of the coagulation equations remain in the self-similar region
$x\approx t^{\frac{1}{1-\gamma}}$ for arbitrarily long times. Since these
estimates are a crucial ingredient in the proof of Theorem \ref{mainRough} we
will describe in this Section how these estimates are derived for the
solutions of the multicomponent coagulation equation (\ref{B4}).

\begin{proposition}
\label{PropMomEst}Suppose that $K:({{\mathbb{R}}_{\ast}^{d}})^{2}%
\rightarrow{\mathbb{R}}_{+}$ is a coagulation kernel satisfying
(\ref{eq:condK_cont}), the homogeneity property (\ref{eq:condK_homogeneity})
as well as the the upper and lower bounds (\ref{eq:condK_sym1}),
\eqref{eq:condK_sym2} with    $\gamma,\gamma+p \in [0,1)$. 
Let $f$ be a solution satisfying   $f(0,\cdot)=f_{0}$ and \eqref{eq:ini_cond} whose existence was stated in Theorem \ref{thm:existence_td}.  Then, for any $k\in{\mathbb{R}}$ satisfying $k\in
(\gamma,1+r]$ if $p\leq0$ or $k\in(-\infty,1+r]$ if $p>0$, there is a constant
$c>0$ that may depend on the initial data such that, for all $t\geq1$,
\begin{equation}
\int_{{\mathbb{R}}_{\ast}^{d}}|x|^{k}f(dx,t)\leq ct^{\frac{k-1}{1-\gamma}}.
\label{MomentEst1}%
\end{equation}

\end{proposition}

The proof of Proposition \ref{PropMomEst} follows directly from the next two lemmas, each
of which provides bounds for the moments $k>1$ and $k<1$ respectively, of a
solution $F$ to the coagulation equation in self-similar variables. More precisely, $F$ satisfies
\begin{align}
\frac{d}{d\tau}\int_{{\mathbb{R}}_{\ast}^{d}}F(d\xi,\tau)\psi(\xi,\tau)  &
=\int_{{{\mathbb{R}}_{\ast}^{d}}}F(d\xi,\tau)\left[  \partial_{\tau}\psi
-\frac{1}{1-\gamma}\xi\cdot\partial_{\xi}\psi+\frac{1}{1-\gamma}\psi\right]
(\xi,\tau)\nonumber\\
&  +\frac{1}{2}\int_{{{\mathbb{R}}_{\ast}^{d}}}\int_{{{\mathbb{R}}_{\ast}^{d}%
}}\tilde{K}\left(  \xi,\eta,\tau\right)  F(d\xi,\tau)F(d\eta,\tau)[\psi
(\xi+\eta,\tau)-\psi(\xi,\tau)-\psi(\eta,\tau)] \label{WeakSelfSimVar}%
\end{align}
for all $\psi\in C_{c}({\mathbb{R}}_{\ast}^{d} \times (0,\infty))$. Notice that
this identity is satisfied for $a.e.$ $\tau\in(0,\infty)$ as it might be seen
using (\ref{eq:weak}) and self-similar variables.

\begin{lemma}
\label{L1}Let $\gamma,p$ and the kernel $K$ satisfy the conditions of
Proposition \ref{PropMomEst}. Let $F_{0}\in\mathcal{M}_{+}({\mathbb{R}}_{\ast
}^{d})$
satisfy
\begin{equation*}
\int_{{\mathbb{R}}_{\ast}^{d}}F_{0}(y)|y|dy=m_{0}. 
\end{equation*}
There is a weak solution $F$ to the coagulation equation in self-similar
variables \eqref{WeakSelfSimVar} with initial condition $F(\cdot,0)=F_{0}$.
Then, for all $k\in(\gamma,1)$ there is a positive constant $w_{k}$ depending
on $k$ such that
\begin{equation*}
\int_{{\mathbb{R}}_{\ast}^{d}}F(d\xi,\tau)\min\{|\xi|,1\}^{k}\leq w_{k}%
,\quad\text{for all }\tau\geq1. 
\end{equation*}

\end{lemma}

\begin{proof}
We generalize the proof of Lemma 3.1 in \cite{EM06} to the multicomponent setting.

We can replace the initial value $F_{0}$ by $F_{\varepsilon,0}$ that is
supported in the region $\{\left\vert x\right\vert \geq\varepsilon\}$ such
that in this region $F_{0,\varepsilon}(dx)=F_{0}(dx)$. Similarly, we can use
an $\varepsilon-$truncated coagulation operator such that each solution
$F_{\varepsilon}$ to \eqref{WeakSelfSimVar} remains supported away from the
origin in $\{\left\vert x\right\vert \geq\varepsilon\}$ for all times $t\geq
0$. All computations that we do in this proof are then fully justified for the
regularized problem. As we will see, the moment estimates derived next will be
uniform in the parameter $\varepsilon$ which allows to conclude their validity
for the original problem taking the limit $\varepsilon\rightarrow0$ at the end
of the argument. Since this argument is standard in the study of coagulation
equations we will not reproduce it here. For simplicity we write in the
following $F$ instead of $F_{\varepsilon}$.

Define the time-independent test functions $\psi(\xi,\tau) = \varphi_{A}(\xi)
= \min(A,|\xi|)^{\ell}$, with $A>0$ and $\ell\in(\gamma,1]$, and
$\tilde\varphi_{A}(\xi,\eta) = \varphi_{A}(\xi+\eta)-\varphi_{A} (\xi
)-\varphi_{A}(\eta)$. Computing $\tilde\varphi_{A}(\xi,\eta) $ yields
\begin{align*}
\tilde\varphi_{A}(\xi,\eta) =
\begin{cases}
|\xi+\eta|^{\ell}-|\xi|^{\ell}-|\eta|^{\ell},\quad\text{for}\quad|\xi|+|\eta|
\leq A\\
A^{\ell}-|\xi|^{\ell}-|\eta|^{\ell},\quad\text{for} \quad|\xi| + |\eta| > A,
\ |\xi|, |\eta| \leq A\\
-|\xi|^{\ell},\quad\text{for} \quad|\xi| \leq A,\ |\eta| > A\\
-|\eta|^{\ell}, \quad\text{for} \quad|y| \leq A,\ |\xi| > A\\
-A^{\ell}, \quad\text{for} \quad|\xi|, |\eta| > A.
\end{cases}
\end{align*}
Note that $\tilde\varphi_{A}(\xi,\eta) \leq0$. Moreover, the following
estimate holds
\begin{equation}
\label{eq:estimate_phitilde}\tilde\varphi_{A}(\xi,\eta) \leq- A^{\ell
}\mathbbm{1}_{\{|\xi|,|\eta| \geq A\}}.
\end{equation}

We also have the following bound for the first terms on the right-hand side of
equation \eqref{WeakSelfSimVar}
\begin{equation}
\label{eq:estimate_derivative}\left[  -\frac{1}{1-\gamma}\xi\cdot\partial
_{\xi}\varphi_{A}+\frac{1}{1-\gamma}\varphi_{A}\right]  (\xi) \leq\frac{
1}{1-\gamma}\varphi_{A}(\xi).
\end{equation}

On the other hand, any kernel in the class considered satisfies the lower
bound
\begin{equation}
\label{eq:boundKernel}K(\xi,\eta) \geq c_{1}(|\xi||\eta|)^{\gamma/2}.
\end{equation}

This follows from the lower bound
\begin{align*}
K(\xi,\eta)  &  \geq c_{1}(|\xi|+|\eta|)^{\gamma}\left(  \frac{|\xi|}%
{|\xi|+|\eta|}\right)  ^{-p}\left(  \frac{|\eta|}{|\xi|+|\eta|}\right)
^{-p}\\
&  =c_{1}(|\xi|+|\eta|)^{\gamma+2p}(|\xi||\eta|)^{-p}=c_{1}(|\xi
||\eta|)^{\gamma/2}(|\xi|+|\eta|)^{\gamma+2p}(|\xi||\eta|)^{-p-\gamma/2}%
\end{align*}
and from the fact that $(|\xi|+|\eta|)^{\gamma+2p}(|\xi||\eta|)^{-p-\gamma
/2}\geq2^{2p+\gamma} \geq1$. To obtain the latter inequality we use that, due
to the Young inequality, $|\xi|^{1/2}|\eta|^{1/2} \leq\frac{1}{2}(|\xi
|+|\eta|)$, as well as the fact that $\gamma+2p\geq0$ (see \eqref{eq:Phi_kernel}).

Using \eqref{eq:estimate_phitilde}, \eqref{eq:estimate_derivative} and
\eqref{eq:boundKernel} it follows from \eqref{WeakSelfSimVar} that, for all
$A>0$,
\begin{equation}
\label{eq:Am}\frac{d}{d\tau}\int_{{\mathbb{R}}_{*}^{d}} F(d\xi,\tau
)\varphi_{A}(\xi) \leq\frac{ 1}{1-\gamma}\int_{{{\mathbb{R}}_{\ast}^{d}}%
}F(d\xi,\tau)\varphi_{A} (\xi) -\frac{A^{\ell}}{2} \left(  \int_{\{|\xi|\geq
A\}} |\xi|^{\gamma/2} F(d\xi,\tau)\right)  ^{2}.
\end{equation}
The strategy now is to obtain a differential inequality for the moment
$\min(1,|\xi|)^{\gamma+\delta}$ with $0< \delta<1-\gamma$.

Define $\phi(y)=\min(1,y)^{\gamma/2+\delta}$. Using integration by parts we
may write
\[
\int_{{\mathbb{R}}_{\ast}^{d}}\phi(|\xi|)|\xi|^{\gamma/2}F(d\xi,\tau)=\int
_{0}^{\infty}\phi^{\prime}(A)\left(  \int_{|\xi|>A}|\xi|^{\gamma/2}F(d\xi
,\tau)\right)  dA
\]
and from Cauchy-Schwarz inequality, we obtain the estimate
\begin{align*}
&  \left(  \int_{0}^{\infty}\phi^{\prime}(A)\left(  \int_{|\xi|>A}
|\xi|^{\gamma/2}F(d\xi,\tau)\right)  dA\right)  ^{2}\leq\\
&  c\int_{0}^{\infty}\phi^{\prime}(A)A^{\gamma/2}\left(  \int_{|\xi|>A}%
|\xi|^{\gamma/2}F(d\xi,\tau)\right)  ^{2}dA,
\end{align*}
with
\[
c:=\int_{0}^{\infty}\phi^{\prime}(A)A^{-\gamma/2}dA=\int_{0}^{1}A^{-1+\delta
}dA<\infty.
\]
Using now \eqref{eq:Am} it follows
\begin{align}
&  \left(  \int_{{\mathbb{R}}_{\ast}^{d}}\phi(|\xi|)|\xi|^{\gamma/2}%
F(d\xi,\tau)\right)  ^{2} \leq c\int_{0}^{\infty}\phi^{\prime}(A)A^{\gamma
/2}\left(  \int_{|\xi|>A}|\xi|^{\gamma/2}F(d\xi,\tau)\right)  ^{2}%
dA\nonumber\\
\leq &  \ 2c\int_{0}^{\infty}\phi^{\prime}(A)A^{\gamma/2}A^{-\ell}\left(
\frac{1}{1-\gamma}\int_{{{\mathbb{R}}_{\ast}^{d}}}F(d\xi,\tau)\varphi_{A}%
(\xi)-\frac{d}{d\tau}\int_{{\mathbb{R}}_{\ast}^{d}}F(d\xi,\tau)\varphi_{A}%
(\xi)\right) \nonumber\\
=  &  \ 2c\left(  \frac{1}{1-\gamma}\int_{{{\mathbb{R}}_{\ast}^{d}}}%
F(d\xi,\tau)\psi(\xi)-\frac{d}{d\tau}\int_{{\mathbb{R}}_{\ast}^{d}}F(d\xi
,\tau)\psi(\xi)\right)  \label{eq:phi_psi}%
\end{align}
with $\psi(\xi)$ defined by
\[
\psi(\xi)=\int_{0}^{\infty}\varphi_{A}(\xi)\phi^{\prime}(A)A^{\gamma/2-\ell
}dA.
\]
Since $\varphi_{A}(\xi)=\min(A,|\xi|)^{\ell}$, then for the choice
$\ell=\gamma+2\delta\leq1,$ one easily concludes that $\psi$ satisfies the
bounds
\begin{equation}
\frac{1}{C}\min(1,|\xi|)^{\gamma+\delta}\leq\psi(\xi)\leq C\min(1,|\xi
|)^{\gamma+\delta} \label{eq:psi_bounds}%
\end{equation}
for some positive constant $C$.
Then \eqref{eq:phi_psi} together with \eqref{eq:psi_bounds} imply an
inequality for the moment $\min(1,|\xi|)^{\gamma+\delta}$,
\begin{align*}
&  \frac{d}{d\tau}\int_{{\mathbb{R}}_{\ast}^{d}}F(d\xi,\tau)\min
(1,|\xi|)^{\gamma+\delta}+\kappa_1\left(  \int_{{\mathbb{R}}_{\ast}^{d}}
F(d\xi,\tau)\min(1,|\xi|)^{\gamma+\delta}\right)  ^{2}\leq\\
&  \kappa_2\int_{{{\mathbb{R}}_{\ast}^{d}}}F(d\xi,\tau)\min(1,|\xi
|)^{\gamma+\delta}%
\end{align*}
for some positive constants $\kappa_1, \kappa_2$.  Integrating this inequality in
time yields the desired uniform in $\varepsilon$ estimate for $F_{\varepsilon
}$ and $\tau\geq1$.
\end{proof}

\begin{lemma}
Let $\gamma, p$ and the kernel $K$ satisfy the conditions of Proposition
\ref{PropMomEst}. Let $F_{0} \in\mathcal{M}_{+}({\mathbb{R}}_{\ast}^{d})$
satisfy
\begin{equation*}
M_{k}:=\int_{{\mathbb{R}}_{*}^{d}} F_{0}(y) |y|^{k} dy <
C
\end{equation*}
and $F$ be the weak solution to the coagulation equation in self-similar
variables \eqref{WeakSelfSimVar} with initial condition $F(\cdot,0)=F_{0}$
obtained in Lemma \ref{L1}. Then, for all $k \in(1,1+\delta)$ there is a
positive constant $w_{k}$ depending on $k$ such that
\begin{equation*}
\sup_{t\geq0}\int_{{\mathbb{R}}_{*}^{d}} F(t,y) |y|^{k}
dy \leq\max\{ w_{k}, M_{k} \}.
\end{equation*}

\end{lemma}

The idea of the proof is to use the test function $\varphi(x)=|x|^{k}$ and to
obtain an estimate in terms of the lower order moments. This idea has been
widely used in the analysis of one component coagulation equation (see for
instance  \cite{EM06} Lemma 3.4 and the book \cite{book}), and it can be immediately adapted to the
multicomponent case. The use of this test function allows to reduce the
estimate for the moment $\int_{{\mathbb{R}}_{\ast}^{d}}F(t,y)|y|^{k}dy$ to the
estimate of moments with an exponent smaller than $k.$ We can then use the
estimate obtained in Lemma \ref{L1}. Since the argument is by now standard, we
will not give more details  here.\

\section{Long time asymptotics for kernels which are constant along any direction}

In this Section we prove Theorem \ref{TheorConstLines}. We need a preliminary
result yielding well-posedness for (\ref{B4}) with kernels satisfying
(\ref{KLines}).

\begin{lemma}
\label{WellPosK}Suppose that the kernel $K$ is as in (\ref{KLines}). Then, for
any $f_{0}\in L^{1}\left(  \mathbb{R}_{\ast}^{d}\right)  $ satisfying
(\ref{F0Mass}) and (\ref{MomHomOne}) there exists a unique solution $f\in
C^{1}\left(  \left(  0,\infty\right)  ;L^{1}\left(  \mathbb{R}_{\ast}%
^{d}\right)  \right)  \cap C\left(  \left[  0,\infty\right)  ;L^{1}\left(
\mathbb{R}_{\ast}^{d}\right)  \right)  $ to (\ref{B4}) in the classical sense
with initial value $f\left(  \cdot,0\right)  =f_{0}\left(  \cdot\right)  .$
The function $f$ is also a weak solution to (\ref{B4}) in the sense of
Definition \ref{def:time-dep_sol}.
\end{lemma}

\begin{proof}
Due to the boundedness of the kernel $K$ we can prove the existence and
uniqueness of a solution $f$ just reformulating (\ref{B4}) as an integral
equation and using a fixed point argument in the space $C\left(  \left[
0,\infty\right)  ;L^{1}\left(  \mathbb{R}_{\ast}^{d}\right)  \right)  .$ The
fact that $f$ is also a weak solution in the sense of Definition
\ref{def:time-dep_sol} follows by multiplying (\ref{B4}) by a test function
$\varphi\left(  x,t\right)  $ and using integration by parts in the variable
$t$ as well as Fubini's Theorem.   These computations are standard,
we refer to the book \cite{book} for further details.
\end{proof}

We now prove Theorem \ref{TheorConstLines}. \smallskip

\begin{proofof}
[Proof of Theorem \ref{TheorConstLines}]For kernels with the form
(\ref{KLines}) and for initial data $f_{0}\in L^{1}\left(  \mathbb{R}_{\ast
}^{d}\right)  $ with the properties stated in Theorem \ref{TheorConstLines},
the conditions in Theorem \ref{mainRoughImproved} are satisfied. Indeed, we
can apply Lemma \ref{WellPosK} and Theorem \ref{thm:existence_td} with
$\gamma=p=0$ with initial data satisfying (\ref{MomHomOne}) to obtain a
solution $f\in C\left(  \left(  0,\infty\right)  ,\mathcal{M}_{+}\left(
\left(  0,\infty\right)  \times\Delta^{d-1}\right)  \right)  $ to (\ref{B4}).
We define $G\in C\left(  \left(  0,\infty\right)  ,\mathcal{M}_{+}\left(
\left(  0,\infty\right)  \times\Delta^{d-1}\right)  \right)  $ by means of
(\ref{MultSelfSimFormT}), (\ref{rhoTheta}), (\ref{PolDist}). Suppose that the
initial data for $G$ is $\bar{G}\in L^{1}\left(  \mathbb{R}_{\ast}^{d}\right)
.$ We will then write $G\left(  \cdot,\cdot,\tau\right)  =S\left(
\tau\right)  \bar{G}\left(  \cdot,\cdot\right)  .$ Notice that
(\ref{MomentEst1}) (or Proposition \ref{PropMomEst}) implies the estimate%
\begin{equation}
\int_{{\mathbb{R}}_{*}\times\Delta^{d-1}}\rho^{k}G(\rho,\theta,\tau
)d\Omega\leq C_{1},\quad k\in\lbrack\frac{1}{a},a]\ \ \text{for some
}a>1\, . \label{MomConst}%
\end{equation}
We recall that $ \tau=\log(t+1)$ with $t\geq 1$,
and again we assume $\tau \geq \ln 2$ throughout.

Let $m_{0}=\left\vert m\right\vert .$ We denote as $\mathcal{N}\left(
\theta_{0};m_{0};C_{1}\right)  $ the subset of $\mathcal{M}_{+}\left(  \left(
0,\infty\right)  \times\Delta^{d-1}\right)  $ that consists in the measures
$\bar{G}$ supported along the line $\left\{  \theta=\theta_{0}\right\}  $ and
satisfying the estimate
\begin{equation}
\int_{{\mathbb{R}}_{*}\times\Delta^{d-1}}\rho^{k}\bar{G}(\rho,\theta)d\Omega\leq C_{1} \label{MomG}%
\end{equation}
(cf.\ (\ref{MomConst})) and having the mass $\int_{{\mathbb{R}}_{*}\times
\Delta^{d-1}}\rho G(\rho,\theta)d\Omega=m_{0}$.   Notice that $\mathcal{N}\left(  \theta_{0}%
;m_{0};C_{1}\right)  $ is a compact subset of $\mathcal{M}_{+}\left(  \left(
0,\infty\right)  \times\Delta^{d-1}\right)  $ in the weak$-\ast$ topology of
$\mathcal{M}_{+}\left(  \left(  0,\infty\right)  \times\Delta^{d-1}\right)  .$
We will denote as $\operatorname*{dist}\left(  \cdot,\cdot\right)  $ a metric
which characterizes the weak$-\ast$ topology of bounded measures in
$\mathcal{M}_{+}\left(  \left(  0,\infty\right)  \times\Delta^{d-1}\right)  .$
We can then apply Theorem \ref{mainRoughImproved} that implies that
\[
\operatorname*{dist}\left(  G(\rho,\theta,\tau),\mathcal{N}\left(  \theta
_{0};m_{0};C_{1}\right)  \right)  \rightarrow0 \text{\quad as \quad}
\tau\rightarrow\infty.
\]

We denote as $G_{0}\left(  \rho;\theta_{0},m_0\right)  $ the measure 
$$G_{0}\left(  \rho;\theta_{0},m_0\right) = \frac
{4}{\left(  Q(\theta_{0})\right)  ^{2}m_{0}}\frac{1}{\rho^{d-1}}\exp\left(
-\frac{2\rho}{Q(\theta_{0})m_{0}}\right).$$ Given $\bar{G}\in\mathcal{N}\left(  \theta_{0};m_{0};C_{1}\right)  $ we can characterize the evolution
semigroup in terms of the corresponding evolution semigroup for the
one-dimensional coagulation evolution. More precisely, we can obtain $S\left(
\tau\right)  \bar{G}$ as the element of $C\left(  \left(  0,\infty\right)  ,
\mathcal{M}_{+}\left(  \left(  0,\infty\right)  \times\Delta^{d-1}\right)
\right)  $ given by $\tau\rightarrow\hat{G}\left(  \left\vert \xi\right\vert
,\tau\right)  \delta\left(  \theta-\theta_{0}\right)  $ where $\hat{G}\left(
\left\vert \xi\right\vert ,\tau\right)  $ is the solution of the one-component
coagulation equation with constant kernel $K=Q\left(  \theta_{0}\right)  $ and initial
value $\bar{G}.$ The existence and uniqueness of $\hat{G}\left(  \left\vert
\xi\right\vert ,\tau\right)  $ follows from \cite{MP04}.\ The results on
\cite{MP04} imply that for any measure $\bar{G}\in\mathcal{N}\left(
\theta_{0};m_{0};C_{1}\right)  $ we have that $S\left(  \tau\right)  \bar
{G}\rightarrow G_{0}\left(  \rho;\theta_{0},m_{0}\right)  $ as $\tau
\rightarrow\infty$ in the weak$-\ast$ topology of $\mathcal{M}_{+}\left(
\left(  0,\infty\right)  \times\Delta^{d-1}\right)  .$ Moreover, the
compactness of $\mathcal{N}\left(  \theta_{0};m_{0};C_{1}\right)  $ implies
that the convergence is uniform. More precisely, for any $\varepsilon>0$ there
exists $T=T\left(  \varepsilon\right)  >0$ such that for any measure $\bar
{G}\in\mathcal{N}\left(  \theta_{0};m_{0};C_{1}\right)  $ we have that
$\operatorname*{dist}\left(  S\left(  \tau\right)  \bar G\left(  \rho
,\theta\right)  ,G_{0}\left(  \rho;\theta_{0},m_{0}\right)  \right)
<\frac{\varepsilon}{2}$ for $\tau\geq T.$ On the other hand, the evolution
equation yields an evolution semigroup that is continuous in the weak$-\ast$
topology of measures with respect to the initial value. We can then argue as
follows in order to prove that the solution $G(\rho,\theta,\tau)$ is at a
distance smaller than $\varepsilon$ from $G_{0}\left(  \rho;\theta_{0}%
,m_{0}\right)  $ for sufficiently large times.

Let $\varepsilon>0$ be an arbitrarily small number. Then, there exists
$T=T\left(  \varepsilon\right)  $ such that%
\begin{equation}
\operatorname*{dist}\left(  S\left(  T\right)  \bar{G},G_{0}\left(
\rho;\theta_{0},m_{0}\right)  \right)  <\frac{\varepsilon}{2} \label{In1}%
\end{equation}
for any $\bar{G}\in\mathcal{N}\left(  \theta_{0};m_{0};C_{1}\right)  .$

On the other hand, the continuity of the semigroup $S\left(  \tau\right)  $
implies that there exists $\delta=\delta\left(  \varepsilon,T\right)  >0$,
that we can assume to satisfy $\delta<\frac{\varepsilon}{2}$ such that, for any  $G_{1}\in\mathcal{M}%
_{+}\left(  \left(  0,\infty\right)  \times\Delta^{d-1}\right)  $ such that
$\operatorname*{dist}\left(  G_{1}(\rho,\theta),\bar{G}\left(  \rho
,\theta\right)  \right)  <\delta$, with $\bar{G}\in\mathcal{N}\left(
\theta_{0};m_{0};C_{1}\right)  $, then
\\ $\operatorname*{dist}\left(  S\left(  T\right)  G_{1}\left(  \rho
;\theta\right)  ,S\left(  T\right)  \bar{G}\left(  \rho,\theta\right)
\right)  <\frac{\varepsilon}{2}.$ Notice that this continuity estimate on the
evolution semigroup is uniform in the class of measures $G_{1}$ satisfying
(\ref{MomG}). The localization result (Theorem \ref{mainRoughImproved})
implies that there exists $T_{1}=T_{1}\left(  \varepsilon,T\right)
=T_{1}\left(  \varepsilon\right)  $ such that for any $\bar{\tau}\geq T_{1}$
we have $\operatorname*{dist}\left(  G(\rho,\theta,\bar{\tau}),\mathcal{N}%
\left(  \theta_{0};m_{0};C_{1}\right)  \right)  <\delta.$ This implies that
there exists $\bar{G}\in\mathcal{N}\left(  \theta_{0};m_{0};C_{1}\right)  $
such that $\operatorname*{dist}\left(  G(\rho,\theta,\bar{\tau}),\bar
{G}\left(  \rho,\theta\right)  \right)  <\delta$  for any
$\bar{\tau}\geq T_{1}.$ Then, given any $\tau\geq T_{1}+T$ we can write
$\tau=\bar{\tau}+T$, which ensures that $\bar{\tau}\geq T_{1}$. Then we obtain 
\begin{align*}
&  \operatorname*{dist}\left(  G\left(  \rho,\theta,\tau\right)  ,G_{0}\left(
\rho;\theta_{0},m_{0}\right)  \right) \\
&  =\operatorname*{dist}\left( S\left(  T\right)  G\left(  \rho
,\theta,\bar{\tau}\right)   ,G_{0}\left(  \rho;\theta_{0},m_{0}\right)  \right) \\
&  \leq\operatorname*{dist}\left(  S\left(  T\right)  G\left(  \rho
,\theta,\bar{\tau}\right)  ,S\left(  T\right)  \bar{G}\left(  \rho,\theta\right)  \right)  +\operatorname*{dist}\left(  S\left(  T\right)
\bar{G}\left(  \rho,\theta\right)  ,G_{0}\left(  \rho;\theta_{0},m_{0}\right)
\right) \\
&  <\frac{\varepsilon}{2}+\frac{\varepsilon}{2}=\varepsilon.
\end{align*}
 Since $\varepsilon$ is arbitrary, the result follows.
\end{proofof}

\begin{remark}
 Combining the methods used in the previous proof with the
ones used in \cite{CT21, Seb2} it would be possible to prove convergence to a self-similar
solution supported along a particular direction for coagulation kernels that
are near constant along each particular direction of the space of cluster compositions.
\end{remark}

\bigskip

\noindent \textbf{Acknowledgements.}
The authors gratefully acknowledge the support of the Hausdorff Research Institute for Mathematics
(Bonn), through the {\it Junior Trimester Program on Kinetic Theory},
of the CRC 1060 {\it The mathematics of emergent effects} at the University of Bonn funded through the German Science Foundation (DFG), 
 of the {\it Atmospheric Mathematics} (AtMath) collaboration of the Faculty of Science of University of Helsinki, of the ERC Advanced Grant 741487 as well as of  
the Academy of Finland via the {\it Centre of Excellence in Analysis and Dynamics Research} (project No.\ 307333).
The funders had no role in study design, analysis, decision to
publish, or preparation of the manuscript.

\noindent\textbf{Compliance with ethical standards} 
\smallskip

\noindent \textbf{Conflict of interest} The authors declare that they have no conflict of interest.

\bigskip 

 \bigskip
 
 \def\adresse{
\begin{description}

\item[M.~A. Ferreira] 
{Department of Mathematics and Statistics, University of Helsinki,\\ P.O. Box 68, FI-00014 Helsingin yliopisto, Finland \\
E-mail:  \texttt{marina.ferreira@helsinki.fi}}

\item[J. Lukkarinen]{ Department of Mathematics and Statistics, University of Helsinki, \\ P.O. Box 68, FI-00014 Helsingin yliopisto, Finland \\
E-mail: \texttt{jani.lukkarinen@helsinki.fi}}

\item[A. Nota:] {Department of Information Engineering, Computer Science and Mathematics,\\ University of L'Aquila, 67100 L'Aquila, Italy \\
E-mail: \texttt{alessia.nota@univaq.it}}

\item[J.~J.~L. Vel\'azquez] { Institute for Applied Mathematics, University of Bonn, \\ Endenicher Allee 60, D-53115 Bonn, Germany\\
E-mail: \texttt{velazquez@iam.uni-bonn.de}}

\end{description}
}

\adresse
 

\begin{thebibliography}{99}                                                                                               


\bibitem {book}J. Banasiak, W. Lamb, P. Lauren\c{c}ot, Analytic methods for
coagulation-fragmentation models, Volume II, CRC Press (2019).


\bibitem{CT21} J.A. Cañizo, S. Throm. The scaling hypothesis for Smoluchowski's coagulation equation with bounded perturbations of the constant kernel. \emph{ J.  Differ. Equ.} \textbf{270} (2021) 285--342.


\bibitem {EM05}M. Escobedo, S. Mischler, M. Rodr\'{\i}guez Ricard, On
self-similarity and stationary problem for fragmentation and coagulation
models, \emph{Ann. Inst. Henri Poincar\`{e} (C) Analyse. Non Lin\'{e}aire},
\textbf{22}(1) (2005) 99--125.

\bibitem {EM06} M. Escobedo, S. Mischler, Dust and self-similarity for the
Smoluchowski coagulation equation. \emph{Ann. Inst. Henri Poincar\`{e} (C)
Analyse. Non Lin\'{e}aire} \textbf{23}(3) (2006) 331--362.

\bibitem {FDGG}J. M. Fern\'{a}ndez-Diaz, G.J. G\'{o}mez-Garc\'{\i}a, Exact
solution of Smoluchowski's continuous multi-component equation with an
additive kernel, \emph{EPL} (Europhysics Letters), \textbf{78}(5) (2007) 56002.

\bibitem {FDGG2}J.M. Fern\'{a}ndez-Diaz, G.J. G\'{o}mez-Garc\'{\i}a. Exact
solution of a coagulation equation with a product kernel in the multicomponent
case. \emph{Phys. D: Nonlinear Phenom.} \textbf{239}(5) (2010) 279--290.


\bibitem{chp}  M.A. Ferreira,  Coagulation Equations for Aerosol Dynamics. In: Albi G., Merino-Aceituno S., Nota A., Zanella M. (eds) \emph{Trails in Kinetic Theory.} SEMA SIMAI Springer Series, vol 25. Springer, Cham. (2021) 69--96.

\bibitem {FLNV} M.A. Ferreira, J. Lukkarinen, A. Nota, J.J.L. Vel\'{a}zquez,
Stationary non-equilibrium solutions for coagulation systems. \emph{Arch.
Ration. Mech. Anal.} \textbf{240} (2021) 809--875 .

\bibitem {FLNV2} M.A. Ferreira, J. Lukkarinen, A. Nota, J.J.L.
Vel\'{a}zquez, Localization in stationary non-equilibrium solutions for
multicomponent coagulation systems. \emph{Comm. Math. Phys.} \textbf{388}(1)
(2021) 479--506.

\bibitem {FLNV4}M.A. Ferreira, J. Lukkarinen, A. Nota, J.J.L. Vel\'{a}zquez,
Multicomponent coagulation systems: existence and non-existence of stationary
non-equilibrium solutions. arXiv:2103.12763 (2021).



\bibitem {FL04}N. Fournier, P. Lauren\c{c}ot, Existence of self-similar solutions
to Smoluchowski's coagulation equation. \emph{Comm. Math. Phys.}
\textbf{256}(3) (2005) 589--609.

\bibitem {FL06}N. Fournier, P. Lauren\c{c}ot. Well-posedness of Smoluchowski's
coagulation equation for a class of homogeneous kernels. \emph{J. of Funct.
Anal.} \textbf{233} (2006) 351--379.

\bibitem {Fried}S.K. Friedlander, Smoke, Dust, and Haze, Oxford University
Press (2000).


\bibitem {KBN}P. Krapivsky, E. Ben-Naim. Aggregation with Multiple
Conservation Laws. \emph{Phys. Rev. E} \textbf{53}(1)  (1995) 291--298.

\bibitem {L} A.A. Lushnikov, Evolution of coagulating systems. III.
Coagulating mixtures, \textit{J. Colloid. Interf. Sci.} \textbf{54} (1) (1976) 94--101.

\bibitem {MP04}G. Menon, R. Pego. Approach to self-similarity in
Smoluchowski's coagulation equation. \emph{Comm. Pure and Appl. Math.}
\textbf{57} (9) (2004) 1197--1232.

\bibitem {MP06}G. Menon, R.L. Pego, Dynamical scaling in Smoluchowski's
coagulation equations: uniform convergence. {\it SIAM Review} \textbf{48}(4) (2006) 745--768.

\bibitem{NTV16} B. Niethammer, S. Throm,  J.J.L. Vel\'{a}zquez. A uniqueness result for self-similar
profiles to Smoluchowski’s coagulation equation revisited. {\it J. Stat. Phys.}, \textbf{164}(2) (2016) 399–-409.


\bibitem {Norris} J.R. Norris, Smoluchowski's coagulation equation:
Uniqueness, nonuniqueness and a hydrodynamic limit for the stochastic
coalescent. \emph{Ann. App. Probab.} \textbf{9}(1) (1999) 78--109.

\bibitem{Olenius} T. Olenius, O. Kupiainen-M\"{a}\"{a}tt\"{a}, I. K. Ortega,
T. Kurt\'{e}n, and H. Vehkam\"{a}ki. Free energy barrier in the growth of sulfuric acid--ammonia and sulfuric acid--dimethylamine clusters. \emph{J. Chem. Phys. }\textbf{139} (2013)  084312.


\bibitem {S16}M. Smoluchowski, Drei vortr\"{a}ge \"{u}ber diffusion, brownsche
molekularbewegung und koagulation von kolloidteilchen. \emph{Phys. Z.} \textbf{17} (1916) 557--585.



\bibitem {Seb2}S. Throm, Stability and Uniqueness of Self-similar Profiles in
$L^{1}$ Spaces for Perturbations of the Constant Kernel in Smoluchowski's
Coagulation Equation. \emph{Commun. Math. Phys.}, \textbf{383}(3) (2021) 1361--1407.


\bibitem {Vehkam} H.~Vehkam\"{a}ki, I.~Riipinen, {Thermodynamics and kinetics
of atmospheric aerosol particle formation and growth\/}, \emph{Chem. Soc.
Rev.}~\textbf{41}(15) (2012) 5160.
\end{thebibliography}
\end{document}